\documentclass[preprint,12pt]{elsarticle}

\usepackage{amssymb}

\usepackage{url}

\usepackage{amssymb}
\usepackage{url}

\usepackage{amsmath,bm,amsfonts}
\usepackage{floatrow}

\usepackage{breqn}

\usepackage{enumerate}
\usepackage{enumitem}

\usepackage{subcaption}

\usepackage{mathtools}

\usepackage{latexsym}

\usepackage{amsthm}
\usepackage{breqn}

\usepackage{enumerate}
\usepackage{enumitem}

\usepackage{bm}

\newtheorem{mydefinition}{Definition}

\def\vE{ {{\bf E}} }
\def\vf{ {{\bf f}} }

\def\vu{ {{\bf u}} }
\def\vv{ {{\bf v}} }

\def\vx{ {{\bf x}} }

\newcommand\scalemath[2]{\scalebox{#1}{\mbox{\ensuremath{\displaystyle #2}}}}

\journal{}

\begin{document}

\begin{frontmatter}

\bibliographystyle{plain}
\bibliographystyle{plainurl}
\nocite{*}


\title{Error Inhibiting Methods for Finite Elements}


\author{Adi Ditkowski}
\address{School of Mathematical Sciences, Tel Aviv University, Tel Aviv 69978, Israel}
\ead{adid@tauex.tau.ac.il}

\author{Anne Le Blanc}

\address{School of Mathematical Sciences, Tel Aviv University, Tel Aviv 69978, Israel}
\ead{anneleb@tauex.tau.ac.il}

\author{Chi-Wang Shu}

\address{Division of Applied Mathematics, Brown University, Providence, RI 02912, USA}
\ead{chi-wang\_shu@brown.edu}

\begin{abstract}
%
%
%

Finite Difference methods (FD) are one of the oldest and simplest methods for solving partial differential equations (PDE).
Block Finite Difference methods (BFD) are FD methods in which the domain is divided into blocks, or cells, containing two or more grid points, with a different scheme used for each grid point, unlike the standard FD method. 

It was shown in \cite{ditkowski2015high} and \cite{ditkowski2020error}
that BFD schemes might be one to three orders more accurate than their truncation errors. Due to these schemes' ability to inhibit the accumulation of truncation errors, these methods were called Error Inhibiting Schemes (EIS). 

This manuscript shows that our BFD schemes can be viewed as a particular type of Discontinuous Galerkin (DG) method. Then, we prove the BFD scheme's stability using the standard DG procedure while using a Fourier-like analysis to establish its optimal convergence rate.

We present numerical examples in one and two dimensions to demonstrate the efficacy of these schemes.


%
%

\end{abstract}


\end{frontmatter}


\section{Introduction}
\label{intro}

This manuscript is concerned with numerical solutions to the Heat equation:
\begin{eqnarray}  \label{Intro_10}
\dfrac{\partial u(\vx,t) }{\partial t} &=& \nabla^2 u(\vx,t) \;, \qquad x\in \Omega \subset  {\mathbb{R}}^d\mbox{ , }t\geq 0  \nonumber\\
u(\vx,t=0)&=& f(\vx)
\end{eqnarray}

Let $Q$ be the discretization of the Laplacian over some discrete, $N$-dimension subspace $\cal{S}$. This subspace could be grid points, as in the case of the Finite Difference (FD), or some finite-dimension function space, as in the Finite elements (FE) or the Discontinuous Galerkin (DG) method. We assume that $Q$ is \textit{semibounded} in the sense that there exists a $H$-norm and a constant $\alpha$ such that 
\begin{equation}  \label{Intro_20}
{\left( Qv,v\right)}_{H}+{\left(v,Qv \right)}_{H}\leq 2\alpha {\parallel v\parallel}^{2}_{H}
\end{equation}

\bigskip

We define the local truncation error of $Q$, $\mathbf{T_{e}}$, by
\begin{equation}  \label{Intro_30}
\left( \mathbf{T_{e}}\right)  := P( \nabla^2 w(x))-(Q\mathbf{w})
\end{equation}
where $w(\vx)$ is a smooth function and $\mathbf{w}$ is the projection of $w(x)$ onto the grid. $P$ is the projection operator onto the subspace $\cal{S}$.\\
 It is assumed that $Q$ is\textit{ consistent } in the sense that 
 \begin{equation}  \label{Intro_40}
\lim_{N\to\infty} {\Vert \mathbf{T_{e}}\Vert}_{H}=0
\end{equation}

Under these conditions, the Lax-Richtmyer equivalence theorem \cite{Richtmyer} assures that the semi-discrete approximation,
\begin{eqnarray}  \label{Intro_50}
\dfrac{\partial \vv(t) }{\partial t} &=&Q \vv(t) \;, \qquad t\geq 0  \nonumber\\
\vv(t=0)&=& \vf(\vx)
\end{eqnarray}
converges. Furthermore, this theorem assures that the error, which is defined as
\begin{eqnarray}  \label{Intro_60}
E = P(u) - \vv
\end{eqnarray}
satisfies 
\begin{eqnarray}  \label{Intro_70}
\left \| E(t) \right \|_H \le K(t) \max_{\tau\in[0,t]} \left \|  \mathbf{T_{e}}(\tau) \right \|_H 
\end{eqnarray}
where $K(\tau) < \text{Const}$ for all  ${\tau\in[0,t]}$, that is,  the error is, at most, of the same order as the truncation error. For most methods, the error is of the same order as the truncation error.

\bigskip

The research published in \citep{ditkowski2015high} and \citep{ditkowski2020error} introduced the Block Finite Difference methods (BFD). BFD are FD methods in which the domain is divided into blocks, or cells, containing two or more grid points, with a different scheme used for each grid point, unlike the standard FD method. It was also shown that for the Heat equation, the BFD schemes might be one or two orders more accurate than their truncation errors. It was also demonstrated that a convergence rate of three orders higher can be achieved using a post-processing filter. Due to these schemes' ability to inhibit the accumulation of truncation errors, these methods were called Error Inhibiting Schemes (EIS).

\bigskip

In this manuscript, we prove that the BFD is a particular type of $p=1$, nodal-based DG scheme. This relationship allows us to establish the stability of this method using the standard approach utilized for DG techniques for both periodic and Dirichlet boundary conditions. Furthermore, we have developed this method to solve the two- and three-dimensional Heat equations.

\bigskip

This paper is constructed as follows. In Section \ref{Sec:BFD}, we present the BFD scheme and prove this BFD scheme is a nodal-based DG method. The stability of this DG method has been proven. This method is then generalized for the Dirichlet boundary conditions and multi-dimensional Heat equation. In \ref{app:Convergence}, we analyze the BFD scheme in the eigenvectors space, find the optimal free parameter, and highlight the benefits of using the post-processing filters. \ref{app:Post-processing} presents the post-processing filters and their implementations.

\section{Block Finite Difference schemes} \label{Sec:BFD}

\noindent  In this section, we present a BFD scheme for the Heat equation with periodic boundary conditions defined as follows:

\begin{equation}
\label{heat_eq_1}
\left\lbrace 
\begin{array}{ll}

\dfrac{\partial u}{\partial t}= \dfrac{\partial ^2}{\partial x^2 } u+F(x,t)\mbox{ , }x\in \left( 0,2\pi \right) \mbox{ , }t\geq 0\\
\\
u(x,0)=f(x)\\

\end{array}
  \right.
\end{equation}

\noindent The one-dimensional scheme was presented in  \cite{ditkowski2020error}. In this manuscript, we give a different proof of stability based on the observation that this BFD scheme is also a particular type of DG method. Later, we will use this observation to derive a multidimensional scheme.

\subsection{Description of Two Points Block, 5th order scheme}\label{Scheme Description}

\noindent Let us consider the following grid:
\begin{equation}
\; x_{j-1/4}=x_{j}-h/4\;,\;x_{j+1/4}=x_{j}+h/4\;,\; h=2\pi/N , \; j=1,...,N
\end{equation}
where 
\begin{equation}
x_{j}=h(j-1)+\frac{h}{2}
\end{equation}

\noindent Altogether there are $2N$ points on the grid, with a distance of $h/2$ between them.
Unlike the standard FD grid, the boundary points do not coincide with any grid nodes.

\noindent We consider the following approximation:
\begin{eqnarray} \label{BFD_scheme_10}
\frac{d^2}{dx^2}u_{j-1/4} &\approx& \frac{1}{12(h/2)^2}\left[ \left( -u_{j-\frac{5}{4}}+16u_{j-\frac{3}{4}}-30u_{j-\frac{1}{4}}+16u_{j+\frac{1}{4}}-u_{j+\frac{3}{4}}\right)+
\right. \nonumber \\
&& \left. c\left( u_{j-\frac{5}{4}} -5u_{j-\frac{3}{4}}+10u_{j-\frac{1}{4}}-10u_{j+\frac{1}{4}}+5u_{j+\frac{3}{4}}-u_{j+\frac{5}{4}}\right)\right] \nonumber \\
\frac{d^2}{dx^2}u_{j+1/4} &\approx& \frac{1}{12(h/2)^2}\left[ \left( -u_{j-\frac{3}{4}}+16u_{j-\frac{1}{4}}-30u_{j+\frac{1}{4}}+16u_{j+\frac{3}{4}}-u_{j+\frac{5}{4}}\right)+
\right. \nonumber \\
&& \left. c\left( -u_{j-\frac{5}{4}} +5u_{j-\frac{3}{4}}-10u_{j-\frac{1}{4}}+10u_{j+\frac{1}{4}}-5u_{j+\frac{3}{4}}+u_{j+\frac{5}{4}}\right)\right]  \nonumber \\
\end{eqnarray}

\noindent or, equivalently:

\begin{eqnarray}\label{2.6}
\vu_{xx} &\approx & \frac{1}{12(h/2)^2}\left[\left(
                                \begin{array}{ccccccccc}
                                     \ddots &\ddots & \ddots & \ddots &  \ddots &   &   &      \\
                                    & -1  & 16 & {\bf -30} & 16 & -1  &   &       \\
                                    &   &  -1 & 16 & {\bf -30} & 16 & -1  &      \\

                                    &   &   & \ddots  & \ddots & \ddots & \ddots  & \ddots  \\
                                \end{array}
                              \right) \right . \\
       &+& c \left . \left(
                                \begin{array}{ccccccccc}
                                  \ddots &   \ddots & \ddots & \ddots &  \ddots &  \ddots &   &      \\
                                   &  1  & -5 & {\bf 10} & -10 & 5  & -1  &       \\
                                   &  -1  & 5 & -10 & {\bf 10} & -5 & 1 &    \\
                                 &      & \ddots  &  \ddots & \ddots & \ddots & \ddots & \ddots   \\
                                \end{array}
                              \right)   \right ] \vu  \; = \;  Q  \vu \nonumber
\end{eqnarray}

\noindent This is a third-order approximation. It was derived using Taylor expansion, with the consistency constraint ${{\bf 1}}^T Q={{\bf 0}}$, where ${{\bf 1}}$ and {{\bf 0}} are column vector with entries 1 and 0 respectively. This consistency constraint  is the numerical equivalence to $\int_0^{2 \pi} u_{xx}dx=0$, for $2 \pi$ periodic functions.


\subsection{Equivalency between the BFD and the  DG Methods}\label{Stability}


%

In this section, we establish the equivalency between the BFD scheme, presented in \eqref{BFD_scheme_10} and \eqref{2.6}, and a particular type of DG method. The properties of the BFD were proven using a Fourier-like analysis in \citep{ditkowski2015high} and \citep{ditkowski2020error}. The equivalency enables us to use the standard DG analysis to prove the method's stability as well. It also enables us to derive high-dimensional efficient DG methods.
 
\subsubsection{The DG Scheme}

Discontinuous Galerkin methods (DG) are a class of finite element methods using discontinuous basis functions \citep{zhang2003analysis}, usually chosen as piecewise polynomials. Consider the following Heat problem,

\begin{equation}
\label{21}
\left\lbrace 
\begin{array}{ll}

\dfrac{\partial u}{\partial t}= \dfrac{\partial ^2}{\partial x^2 } u\mbox{ , }x\in \left( 0,2\pi \right) \mbox{ , }t\geq 0\\
\\
u(x,0)=f(x)\\

\end{array}
  \right.
\end{equation}

We assume the following mesh to cover the computational domain $[0,2\pi]$, consisting of cells

$I_{j}=\left[ x_{j-{1/2}},x_{j+{1/2}}\right] $ for $j=1,...,N$, where\\
\begin{center}

$0=x_{1/2}<x_{3/2}<...<x_{N+1/2}=2\pi.$\\
\end{center}

\noindent The center of each cell is located at $x_{j}=\frac{1}{2}\left( x_{j-{1/2}}+x_{j+{1/2}}\right) $ and the size of each cell is $\Delta x_{j}=x_{j+{1/2}}-x_{j-{1/2}}$. A uniform mesh is considered here, hence $h=\Delta x=\frac{2\pi}{N}$.\\

\noindent After multiplying the two sides of the equation by a test function, $v(x)$,  and integrating by parts over each cell $I_{j}$, we get the following weak formulation of the problem:
\begin{equation}
\label{2}
\left.
\begin{array}{ll}
\int_{I_j}u_{t}v dx+\int_{I_j}u_{x}v_{x}dx-u_{x}(x_{j+1/2},t)v(x_{j+1/2})+u_{x}(x_{j-1/2},t)v(x_{j-1/2})=0
  \end{array}
  \right.
\end{equation}

\noindent The next step of the method is to replace the functions $u$ and $v$ at each cell by piecewise polynomials of degree at most $k$. We still denote $u$ and $v$ for the approximated functions to simplify the notation. \\
\noindent Since both functions are now discontinuous at all points $x_{j\pm 1/2}$, there is a need to introduce a value for $u_{x}(x_{j\pm 1/2},t)$ and $v(x_{j\pm 1/2},t)$.  We replace the boundary terms $u_{x}(x_{j\pm 1/2},t)$ by single-valued fluxes $\hat{u}_{j\pm 1/2}=\hat{u}((u_x)_{j\pm 1/2}^{-},(u_x)_{j\pm 1/2}^{+})$ and replace the test function $v$ at the boundaries by its values taken inside the cell. $\hat{u}_{j\pm 1/2}$ are called \textit{numerical fluxes}. Here, unlike the hyperbolic case, there is no preferred flux direction. Therefore, we choose a central flux:  
\begin{equation}
\label{3}
\left.
\begin{array}{ll}
\hat{u}_{x,j+\frac{1}{2}}=\frac{1}{2}\left[ (u_{x})^{-}_{j+\frac{1}{2}}+(u_{x})^{+}_{j+\frac{1}{2}}\right] \\
\hat{u}_{x,j-\frac{1}{2}}=\frac{1}{2}\left[ (u_{x})^{-}_{j-\frac{1}{2}}+(u_{x})^{+}_{j-\frac{1}{2}}\right]
  \end{array}
  \right.
\end{equation}

\noindent The corresponding DG scheme is 
\begin{equation}
\label{2}
\left.
\begin{array}{ll}
\int_{I_j}u_{t}v dx+
		\int_{I_j}u_{x} v_{x}dx- \\
\hspace {1.5em}		\frac{1}{2}\left[ (u_{x})^{-}_{j+\frac{1}{2}}+(u_{x})^{+}_{j+\frac{1}{2}}\right]v^{-}(x_{j+1/2})+\frac{1}{2}\left[ (u_{x})^{-}_{j-\frac{1}{2}}+(u_{x})^{+}_{j-\frac{1}{2}}\right]v^{+}(x_{j-1/2})=0
  \end{array}
  \right.
\end{equation}
\noindent Unfortunately, this scheme is unstable \citep{zhang2003analysis}. In order to make it so, it is necessary to add penalty terms to inter-element boundaries, namely the  Baumann-Oden penalty terms (see \citep{Baumann_Oden} for more details). We note here that another possible modification can be made to achieve stability (see \cite{zhang2003analysis}).

The DG scheme is then defined as follows:

\noindent Find $u\in V_{\Delta x}$ (where $V_{\Delta x}=\lbrace v:$ $v$ is a polynomial of degree at most $k$ for $x\in I_{j},j=1,...,N\rbrace$) such that for all $v\in V_{\Delta x}$,
\begin{equation}
\label{4}
\left.
\begin{array}{ll}
\int_{I_j}u_{t}v dx+\int_{I_j}u_{x}v_{x}dx-\hat{u}_{x}(x_{j+1/2},t)v^{-}(x_{j+1/2})+\hat{u}_{x}(x_{j-1/2},t)v^{+}(x_{j-1/2})\\
\hspace {3em}		
\underbrace{ -\frac{1}{2}v_{x}^{-}(x_{j+1/2})\left[u^{+}_{j+1/2}-u^{-}_{j+1/2} \right]-
\frac{1}{2}v_{x}^{+}(x_{j-1/2})\left[u^{+}_{j-1/2}-u^{-}_{j-1/2} \right] }_{\text{Baumann-Oden penalty term}  }= 0
  \end{array}
  \right.
\end{equation}


Following \cite{zhang2003analysis}, we now look at the nodal presentation of the scheme; namely, we express the functions $u(x,t)$ and $v(x)$ at the nodes $(x_{j \pm 1/4})
$. When a linear element basis is used on an equidistant grid, $\varphi_{j-1/4}$ and $\varphi_{j+1/4}$, $u$ and $v$ have the form
\begin{eqnarray} \label{nodal_u_v_def}
u &=& u_{j-1/4} \varphi_{j-1/4}  + u_{j+1/4} \varphi_{j+1/4} \nonumber \\
v &=& v_{j-1/4} \varphi_{j-1/4}  + v_{j+1/4} \varphi_{j+1/4}
\end{eqnarray}
\noindent where
\begin{equation}
\label{66}
\left.
\begin{array}{ll}

   \varphi_{j-1/4}= -\frac{2}{h}(x-x_{j+\frac{1}{4}})\\
    \varphi_{j+1/4}=\frac{2}{h}(x-x_{j-\frac{1}{4}})\\

  \end{array}
  \right.
\end{equation}
\noindent are the Lagrange interpolating polynomials. Then,\\
\begin{equation}\nonumber
\label{9_new}
\left.
\begin{array}{ll}
(u)^{+}_{j+1/2}=u_{j+3/4}\varphi_{j+3/4}(x_{j+1/2})+u_{j+5/4}\varphi_{j+5/4}(x_{j+1/2})\\
(u)^{-}_{j+1/2}=u_{j-1/4}\varphi_{j-1/4}(x_{j+1/2})+u_{j+1/4}\varphi_{j+1/4}(x_{j+1/2})\\
(u)^{+}_{j-1/2}=u_{j-1/4}\varphi_{j-1/4}(x_{j-1/2})+u_{j+1/4}\varphi_{j+1/4}(x_{j-1/2})\\
(u)^{-}_{j-1/2}=u_{j-5/4}\varphi_{j-5/4}(x_{j-1/2})+u_{j-3/4}\varphi_{j-3/4}(x_{j-1/2})\\
  \end{array}
  \right.
\end{equation}
\noindent Collecting the coefficients of $v_{j-1/4}$ and $v_{j+1/4}$ yields the equation for $u_{j-1/4}$ and $u_{j+1/4}$
\begin{equation}
\label{5}
\left.
\begin{array}{ll}

   \begin{bmatrix}
           u_{j-1/4} \\
           u_{j+1/4} \\
          
         \end{bmatrix}_{t} &=\left( A\begin{bmatrix}
           u_{j-5/4} \\
           u_{j-3/4} \\
          
         \end{bmatrix} +B\begin{bmatrix}
           u_{j-1/4} \\
           u_{j+1/4} \\
          
         \end{bmatrix}+C\begin{bmatrix}
           u_{j+3/4} \\
           u_{j+5/4} \\
          
         \end{bmatrix}
 \right)
  \end{array}
  \right.
\end{equation}
\noindent where
\begin{eqnarray} \label{6_2.5}
 A &=& \frac{1}{4 h^2}\begin{bmatrix}
          7 & -1\\
          1 & -7 \\
         \end{bmatrix}   \nonumber \\
B &=& \frac{1}{2 h^2}\begin{bmatrix}
           -12 & 12 \\
          12 & -12\\
        \end{bmatrix}         \\
C &=& \frac{1}{4 h^2}\begin{bmatrix}
          -7 & 1\\
          -1 & 7 \\
           \end{bmatrix}         \nonumber
\end{eqnarray}
\noindent This scheme has been proven to be consistent, stable, and of $k$ order accuracy for even $k$ and $k+1$ for odd $k$ (\cite{zhang2003analysis} \cite{BABUSKA1999103}).\\

\noindent Clearly, we can write the BFD scheme as defined in Section \ref{Scheme Description} in a similar form by assembling the matrices $A$, $B$, and $C$  as follows:
\begin{eqnarray} \label{6_3}
 A &=& \frac{1}{3h^2}\begin{bmatrix}
           -1+c&16-5c \\
           -c&-1+5c \\
         \end{bmatrix}   \nonumber \\
B &=& \frac{1}{3h^2}\begin{bmatrix}
           -30+10c&16-10c \\
           16-10c&-30+10c\\
        \end{bmatrix}         \\
C &=& \frac{1}{3h^2}\begin{bmatrix}
           -1+5c&-c \\
           16-5c&-1+c \\
           \end{bmatrix}         \nonumber
\end{eqnarray}
\noindent Hence, our goal is to find the problem's corresponding weak formulation, including the Baumann-Oden and other penalty terms, as well as numerical fluxes, such that our BFD scheme can be viewed as a form of DG scheme.

\noindent As done in \cite{zhang2003analysis}, we choose a linear element basis \eqref{nodal_u_v_def} for the test and trial functions. By replacing in Eq. (\ref{4}) the fluxes and standard Baumann-Oden penalties by fluxes and all the possible penalties with general coefficients, we obtain the following scheme:

\begin{equation}
\label{8}
\left.
\begin{array}{ll}
\int\limits_{x_{j-1/2}}^{x_{j+1/2}}\left[ (u_{j-1/4})_{t}\varphi_{j-1/4}+(u_{j+1/4})_{t}\varphi_{j+1/4}\right] v(x)dx=\\
 \hspace{2em} -\int\limits_{x_{j-1/2}}^{x_{j+1/2}}\left[ u_{j-1/4}(\varphi_{j-1/4})_{x}+u_{j+1/4}(\varphi_{j+1/4})_{x}\right] v_{x}(x)dx \; +\\
 \hspace{4em} \hat{u}_{x,j+1/2}v^{-}(x_{j+1/2})-\hat{u}_{x,j-1/2}v^{+}(x_{j-1/2})  \; +\\
\\
\Bigg( \dfrac{C_{1} }{h}\bigg((u)^{+}_{j+1/2}-(u)^{-}_{j+1/2}\bigg)+{C}_{2}\bigg((u_{x})^{+}_{j+1/2}-(u_{x})^{-}_{j+1/2}\bigg)\Bigg)v^{-}_{j+1/2}\ \; - \\ \\
\Bigg(  \dfrac{D_{1} }{h}\bigg((u)^{+}_{j-1/2}-(u)^{-}_{j-1/2}\bigg)+{D}_{2}\bigg((u_{x})^{+}_{j-1/2}-(u_{x})^{-}_{j-1/2}\bigg)\Bigg)v^{+}_{j-1/2} \; + \\
\\
\Bigg(E_{1}\bigg((u)^{+}_{j+1/2}-(u)^{-}_{j+1/2}\bigg)+h E_{2}\bigg((u_{x})^{+}_{j+1/2}-(u_{x})^{-}_{j+1/2}\bigg) \Bigg)(v_{x})_{j+1/2}^{-} \; -\\
\\
\Bigg( F_{1}\bigg((u)^{+}_{j-1/2}-(u)^{-}_{j-1/2}\bigg)+h F_{2}\bigg((u_{x})^{+}_{j-1/2}-(u_{x})^{-}_{j-1/2}\bigg)\Bigg)(v_{x})_{j-1/2}^{+}\\
  \end{array} 
  \right.
\end{equation}

\noindent where the following general fluxes were defined by
\begin{equation}
\label{3_1}
\left.
\begin{array}{ll}
\hat{u}_{x,j+1/2}=\alpha(u_{x})^{-}_{j+1/2}+\left( 1-\alpha\right) (u_{x})^{+}_{j+1/2} \\
\\
\hat{u}_{x,j-1/2}=\beta(u_{x})^{-}_{j-1/2}+\left( 1-\beta\right) (u_{x})^{+}_{j-1/2}\\
  \end{array}
  \right.
\end{equation}
\noindent for some  $\alpha,\beta \in \mathbb{R}$. The design of the above scheme is based on the minimal requirement that for a smooth derivable solution, all penalties should cancel each other. By comparing the coefficients of $u_{j-5/4}, ..., u_{j+5/4}$ to the BFD scheme we obtain:
%
%
%
%
%
\begin{equation}
\label{coeff_solution_10}
\left.
\begin{array}{llllll}
C_1=\frac{7}{3}, &&  C_2=\alpha-\frac{1}{2},  \\
\\
D_1=\frac{7}{3}, &&  D_2=\beta-\frac{1}{2}, \\
 \\
 E_1=-\frac{1}{18} (8 c+5),  &&  E_2=- \frac{1}{18} (c+1), \\
 \\
 F_1= \frac{1}{18} (8 c+5), &&  F_2= -\frac{1}{18} (c+1)
 
 \end{array} 
  \right.
\end{equation}
This leaves us with two remaining unknowns to be identified, namely $\alpha$ and $\beta$. In order to determine those, we use the stability tool developed in \citep{shu2009discontinuous} to prove the stability of FE methods. First, we define the following operator:

%

\begin{mydefinition}[The Operator ${\Theta}_{j-1/2}$]

The following operator includes the net contribution viewed as the difference between penalties from each side of the left border of the cell $I_{j}$, the numerical flux applied to the same node at node $x_{j-1/2}$, and the contribution from the integration over both cells when we set $v=u$.


\begin{equation}
\label{3_2}
\left.
\begin{array}{ll}
{\Theta}_{j-1/2}=\bigg(\alpha(u_{x})^{-}_{j-1/2}+\left( 1-\alpha\right) (u_{x})^{+}_{j-1/2} \bigg)u^{-}(x_{j-1/2})\\
\\
-\bigg(\beta(u_{x})^{-}_{j-1/2}+\left( 1-\beta\right) (u_{x})^{+}_{j-1/2}\bigg)u^{+}(x_{j-1/2})\\
\\
+\Bigg( \dfrac{C_{1}}{h} \bigg((u)^{+}_{j-1/2}-(u)^{-}_{j-1/2}\bigg)+{C}_{2}\bigg((u_{x})^{+}_{j-1/2}-(u_{x})^{-}_{j-1/2}\bigg)\Bigg)u^{-}_{j-1/2}\\
\\
-\Bigg(  \dfrac{D_{1}}{h}  \bigg((u)^{+}_{j-1/2}-(u)^{-}_{j-1/2}\bigg)+{D}_{2}\bigg((u_{x})^{+}_{j-1/2}-(u_{x})^{-}_{j-1/2}\bigg)\Bigg)u^{+}_{j-1/2}\\
\\
+\Bigg(E_{1}\bigg((u)^{+}_{j-1/2}-(u)^{-}_{j-1/2}\bigg)+h E_{2}\bigg((u_{x})^{+}_{j-1/2}-(u_{x})^{-}_{j-1/2}\bigg) \Bigg)(u_{x})_{j-1/2}^{-}\\
\\
-\Bigg( F_{1}\bigg((u)^{+}_{j-1/2}-(u)^{-}_{j-1/2}\bigg)+h F_{2}\bigg((u_{x})^{+}_{j-1/2}-(u_{x})^{-}_{j-1/2}\bigg)\Bigg)(u_{x})_{j-1/2}^{+}\\
\\
-\dfrac{1}{2}\int_{I_{j-1}}(u_{x})^2_{j-1}dx-\dfrac{1}{2}\int_{I_{j}}(u_{x})^2_{j}dx\nonumber
  \end{array}
  \right.
\end{equation}

\end{mydefinition}

\noindent The condition that the quadratic form ${\Theta}_{j-1/2}$ is non-positive definite is sufficient for stability.
%

\noindent The operator $\Theta_{j-1/2}$ may be written in the following form : 

\begin{equation}
\label{3_2_1}
\left.
\begin{array}{ll}
{\Theta}_{j-1/2}=\sum_{k,l=1}^{4}a_{k,l}x_{k}x_{l}\nonumber
  \end{array}
  \right.
\end{equation}
\noindent where
\begin{equation}
\label{3_2_2}
\left.
\begin{array}{ll}
\mathbf{x}=\begin{pmatrix} u_{j-5/4} \\  
   u_{j-3/4}\\
     u_{j-1/4}\\
      u_{j+1/4}\end{pmatrix}\nonumber
  \end{array}
  \right.
\end{equation}
\noindent The matrix for the associated symmetric bilinear form
\begin{equation} \label{Theta_expansion}
 {\Theta}_{j-1/2}= 
    \mathbf{x}^{T} 
 M
  \mathbf{x}\nonumber
\end{equation}
is 
\[
 M= \frac{1}{12h}\left[
    \scalemath{1}{
    \begin{array}{cccccccc}
     8 c-19 & \frac{1}{3} (71-40 c) & 8 c-5 & \frac{1}{3} (1-8 c) \\
     \\
    \frac{1}{3} (71-40 c) & \frac{1}{3} (56 c-169) & \frac{1}{3} (113-40 c) & 8 c-5 \\
    \\
     8 c-5 & \frac{1}{3} (113-40 c) & \frac{1}{3} (56 c-169) & \frac{1}{3} (71-40 c) \\
     \\
     \frac{1}{3} (1-8 c) & 8 c-5 & \frac{1}{3} (71-40 c) & 8 c-19 \\
     \\
    \end{array}
    }
  \right]\nonumber
\]

\noindent We note here that $M$ does not depend either on $\alpha$ or $\beta$; hence those are free parameters.

\noindent $M$ is a singular matrix, since $\Theta_{j-1/2}(u)=0$ for any constant function $u$. 

\noindent After checking the rows and columns' linear dependency, the fourth row and columns depend on the other rows and columns, respectively.

\noindent Therefore, it is sufficient to ensure the non-positiveness of the truncated matrix $M'$:
\[M'= \frac{1}{12h}
\left[
    \scalemath{1}{
    \begin{array}{cccccccc}
     
 8 c-19 & \frac{1}{3} (71-40 c) & 8 c-5 \\
 \\
 \frac{1}{3} (71-40 c) & \frac{1}{3} (56 c-169) & \frac{1}{3} (113-40 c) \\
 \\
  8 c-5 & \frac{1}{3} (113-40 c) & \frac{1}{3} (56 c-169) \\

    \end{array}
    }
  \right]
\]
\noindent Using Sylvester's law of inertia \cite{sylvester}, stating that the number of eigenvalues of each sign is an invariant of the associated quadratic form, we use a congruence relation  with $M'$, thus obtaining an invertible matrix $S$ such that $D = SM'S^{T}$, where $D$ is a diagonal matrix.
\noindent Equivalently:
\[D= \frac{1}{12h}
 \left[
    \scalemath{0.8}{
    \begin{array}{cccccccc}
    8 c-19 & 0 & 0  \\
    \\
    0 &  \hspace{-0em} \frac{1}{9} (-16) (8 c-19) (2 c (8 c+49)-287) & 0 \\
    \\
    0 & 0 &  \hspace{-2em} \begin{array}{l}
   		   \frac{16384}{243} (19-8 c)^4 (2 c-7) (2 c (2 c+7)-35) \cdot \\
   		   		 \hspace{2em}  (2 c (8 c+49)-287)  
   		   \end{array} \\
    \end{array}
    }
  \right]
\]

\noindent whose diagonal elements are non-positive for all $-1\leq c \leq 1$, as required.

\bigskip

In the previous section, we proved that our Two Points Block 5th order scheme can be viewed as a type of DG scheme. We also showed that the stability of this scheme can be proven using the tools developed for DG methods.

In the following section, we generalize the principles developed above for a scheme approximating non-periodic problems.

\subsection{Generalization to non periodic boundary conditions}


Consider the following Heat Initial Boundary Value Problem with Dirichlet boundary conditions:
\begin{equation}
\label{21_dir}
\left\lbrace 
\begin{array}{ll}

\dfrac{\partial u}{\partial t}= \dfrac{\partial ^2}{\partial x^2 } u+F(x,t)\mbox{ , }x\in \left( 0,\pi \right) \mbox{ , }t\geq 0\\
\\
u(x,t=0)=f(x)\\
u(0,t)=g_{0}(t)\\
u(\pi,t)=g_{\pi}(t)
\end{array}
  \right.
\end{equation}

\subsubsection{Adapting the BFD scheme to non-periodic boundary conditions of Dirichlet type}

\noindent As the problem \eqref{21_dir}  contains boundary conditions, an adaptation from our fifth-order BFD scheme with periodic conditions is required. As shown in \citep{ditkowski2020error}, the scheme can be applied on the interval $\left( 0,\pi\right)$, then an anti-symmetric reflection is performed onto the interval $\left( 0,2\pi\right)$.\\

\noindent In order to get an approximation scheme for the end points of the interval $[0,\pi]$, we perform an extrapolation to the two additional ghost points needed, $x_{\frac{1}{4}}=-\frac{h}{4}$, $x_{-\frac{1}{4}}=-\frac{3h}{4}$ at the left boundary and $x_{(N+1)-\frac{1}{4}}=\pi+\frac{h}{4}$, $x_{(N+1)+\frac{1}{4}}=\pi+\frac{3h}{4}$ at the right one. As for the internal points for $2\leq j \leq N-1$, the scheme remains the same.

\noindent The extrapolations, using Taylor's expansions, are:
\begin{equation}
\label{dir1}
\left.
\begin{array}{ll}

u_{-\frac{1}{4}}=-u_{1+\frac{1}{4}}+2g_{0}+u_{xx}\left( 0,t\right)\left( \frac{3h}{4}\right)^2+\frac{1}{12}u_{xxxx}\left( 0,t\right)\left( \frac{3h}{4}\right)^4+O(h^6)\\

u_{\frac{1}{4}}=-u_{1-\frac{1}{4}}+2g_{0}+u_{xx}\left( 0,t\right)\left( \frac{h}{4}\right)^2+\frac{1}{12}u_{xxxx}\left( 0,t\right)\left( \frac{h}{4}\right)^4+O(h^6)\\

\end{array}
  \right.
\end{equation}
\noindent and
\begin{equation}
\label{dir2.1}
\left.
\begin{array}{ll}

u_{(N+1)-\frac{1}{4}}=-u_{N+\frac{1}{4}}+2g_{\pi}+u_{xx}\left( \pi,t\right)\left( \frac{h}{4}\right)^2+\frac{1}{12}u_{xxxx}\left( \pi,t\right)\left( \frac{h}{4}\right)^4+O(h^6)\\

u_{(N+1)+\frac{1}{4}}=-u_{N-\frac{1}{4}}+2g_{\pi}+u_{xx}\left( \pi,t\right)\left( \frac{3h}{4}\right)^2+\frac{1}{12}u_{xxxx}\left( \pi,t\right)\left( \frac{3h}{4}\right)^4+O(h^6)\\

\end{array}
  \right.
\end{equation}

\noindent where $u_{xx}(0,t)$,$u_{xxxx}(0,t)$,$u_{xx}(\pi,t)$,$u_{xxxx}(\pi,t)$ can be computed from the PDE.
\begin{equation}
\label{dir_eq1}
\left.
\begin{array}{ll}

u_{xx}(0,t)=u_{t}(0,t)-F(0,t)\\
u_{xx}(\pi,t)=u_{t}(\pi,t)-F(\pi,t)\\
u_{xxxx}(0,t)=u_{tt}(0,t)-F_{t}(0,t)-F_{xx}(0,t)\\
u_{xxxx}(\pi,t)=u_{tt}(\pi,t)-F_{t}(\pi,t)-F_{xx}(\pi,t)\\

\end{array}
  \right.
\end{equation}
Since $u$ satisfies the prescribed boundary conditions $u(0,t)=g_{0}(t)$ and $u(\pi,t)=g_{\pi}(t)$, the time derivatives can be computed analytically or numerically. 

\noindent The schemes at the nodes $x_{1-1/4}$ and $x_{1+1/4}$ become:
\begin{eqnarray}
\label{dir3}
&& \frac{d^2}{dx^2}u_{1-\frac{1}{4}}  \approx \frac{1}{12(h/2)^2}\left[ \left( 30-8c\right)g_{0}+\left( 7+4c\right)\left( \frac{h}{4}\right)^2 u_{xx}(0,t)+ \right . \nonumber \\
&& \hspace{2 em}  \frac{\left( -65+76c\right)}{12} \left( \frac{h}{4}\right)^4 u_{xxxx}(0,t)+ \left(-46u_{1-\frac{1}{4}}+17 u_{1+\frac{1}{4}}-u_{1+\frac{3}{4}}\right)+ 
\nonumber \\
&& \hspace{3 em} \left . c\left( 15u_{1-\frac{1}{4}}-11u_{1+\frac{1}{4}}+5u_{1+\frac{3}{4}}-u_{1+\frac{5}{4}}\right)\right]
\end{eqnarray}

\begin{eqnarray}
\label{dir4}
&& \frac{d^2}{dx^2}u_{1+\frac{1}{4}}\approx\frac{1}{12(h/2)^2}\left[ \left( -2+8c\right)g_{0}+\left(-1-4c\right)\left( \frac{h}{4}\right)^2 u_{xx}(0,t)+ \right . \nonumber \\
&& \hspace{2 em}  \frac{\left( -1-76c\right)}{12} \left( \frac{h}{4}\right)^4 u_{xxxx}(0,t)+ \left(17u_{1-\frac{1}{4}}-30u_{1+\frac{1}{4}}+16u_{1+\frac{3}{4}}-u_{1+\frac{5}{4}}\right)+
\nonumber \\
&& \hspace{3 em} \left .  c\left(-15u_{1-\frac{1}{4}}+11u_{1+\frac{1}{4}}-5u_{1+\frac{3}{4}}+u_{1+\frac{5}{4}}\right)\right]
\end{eqnarray}

\noindent Similarly, the schemes at the nodes $x_{N-1/4}$ and $x_{N+1/4}$ are
\begin{eqnarray}
\label{dir5}
&& \frac{d^2}{dx^2}u_{N-\frac{1}{4}}\approx\frac{1}{12(h/2)^2}\left[ \left( -2+8c\right)g_{\pi}+\left( -1-4c\right)\left( \frac{h}{4}\right)^2 u_{xx}(\pi,t)+
\right . \nonumber \\
&& \hspace{2 em} \frac{\left( -1-76c\right)}{12} \left( \frac{h}{4}\right)^4 u_{xxxx}(\pi,t)+
\left(-u_{N-\frac{5}{4}}+16 u_{N-\frac{3}{4}}-30u_{N-\frac{1}{4}}+17u_{N+\frac{1}{4}}\right)+
\nonumber \\
&& \hspace{3 em} \left .  c\left( u_{N-\frac{5}{4}}-5u_{N-\frac{3}{4}}+11u_{N-\frac{1}{4}}-15u_{N+\frac{1}{4}}\right)\right]
\end{eqnarray}

\begin{eqnarray}
\label{dir5}
&& \frac{d^2}{dx^2}u_{N+\frac{1}{4}}  \approx \frac{1}{12(h/2)^2}\left[ \left( 30-8c\right)g_{\pi}+\left( 7+4c\right)\left( \frac{h}{4}\right)^2 u_{xx}(\pi,t)+ \right . \nonumber \\
&& \hspace{2 em}  \frac{\left( -65+76c\right)}{12} \left( \frac{h}{4}\right)^4 u_{xxxx}(\pi,t)+ \left(-u_{N-\frac{3}{4}}+17 u_{N-\frac{1}{4}}-46u_{N+\frac{1}{4}}\right)+ 
\nonumber \\
&& \hspace{3 em} \left . c\left( -u_{N-\frac{5}{4}}+5u_{N-\frac{3}{4}}-11u_{N-\frac{1}{4}}+15u_{N+\frac{1}{4}}\right)\right]
\end{eqnarray}

where $u_{xx}(0,t),u_{xxxx}(0,t),u_{xx}(\pi,t),u_{xxxx}(\pi,t)$ are computed by \eqref{dir_eq1}.
\\

\noindent The truncation errors at the boundaries are as follows:
\begin{dmath*}
T_e\left( x_{ 1-\frac{1}{4} } \right) =- \frac{2359}{4423680}h^4u^{(6)}+c\left[ -\frac{h^3}{96}u^{(5)}-\frac{3061}{1105920}h^4u^{(6)}\right] +O(h^5)=O(h^3)\\
\end{dmath*}

\begin{dmath*}
T_e\left( x_{1+\frac{1}{4} } \right)=- \frac{3071}{4423680}h^4u^{(6)}+c\left[ \frac{h^3}{96}u^{(5)}-\frac{2699}{1105920}h^4u^{(6)}\right] +O(h^5)=O(h^3)\\
\end{dmath*}

\begin{dmath*}
T_e\left( x_{N-\frac{1}{4} } \right)=- \frac{3071}{4423680}h^4u^{(6)}+c\left[ -\frac{h^3}{96}u^{(5)}-\frac{2699}{1105920}h^4u^{(6)}\right] +O(h^5)=O(h^3)\\
\end{dmath*}

\begin{dmath*}
T_e\left( _{N+\frac{1}{4} } \right)=- \frac{2359}{4423680}h^4u^{(6)}+c\left[ \frac{h^3}{96}u^{(5)}-\frac{3061}{1105920}h^4u^{(6)}\right] +O(h^5)=O(h^3)\\
\end{dmath*}

\subsubsection{Equivalence of the BFD and DG schemes}

Again, stability is proven using the equivalence between our modified BFD scheme and a specific non-standard DG scheme.\\
Now, like in the case of periodic BC, our goal is to find the corresponding weak formulation of the problem, including the Baumann-Oden and possibly other penalty terms as well as numerical fluxes, so that our BFD scheme can be viewed as a form of DG scheme.\\

\noindent In the first cell of the grid $[0,h]$, the scheme can be written as:

\begin{equation}
\label{}
\left.
\begin{array}{ll}

   \begin{bmatrix}
           u_{1-1/4} \\
           u_{1+1/4} \\
          
         \end{bmatrix}_{t} &=\left( B\begin{bmatrix}
           u_{1-1/4} \\
           u_{1+1/4} \\
          
         \end{bmatrix} +C\begin{bmatrix}
           u_{1+3/4} \\
           u_{1+5/4} \\
        
         \end{bmatrix}
 \right)
  \end{array}
  \right.
\end{equation}
where:
\begin{equation}
\label{Chap3_1}
\left.
\begin{array}{ll}

   B= \frac{1}{3h^2}\begin{bmatrix}
           -46+15c&17-11c \\
           17-15c&-30+11c \\
          
         \end{bmatrix} \mbox{ and }
         C= \frac{1}{3h^2}\begin{bmatrix}
           -1+5c&-c \\
           16-5c&-1+c\\
          
         \end{bmatrix}

  \end{array}
  \right.
\end{equation}
\noindent The general scheme would be of the form:
\begin{equation}
\label{Chap3_2}
\left.
\begin{array}{ll}
\int\limits_{x_{1-1/2}}^{x_{1+1/2}}\left[ (u_{1-1/4})_{t}\varphi_{1-1/4}+(u_{1+1/4})_{t}\varphi_{1+1/4}\right] v(x)dx \; =\\
\hspace{2em}
-\int\limits_{x_{1-1/2}}^{x_{1+1/2}}\left[ u_{1-1/4}(\varphi_{1-1/4})_{x}+u_{1+1/4}(\varphi_{1+1/4})_{x}\right] v_{x}(x)dx \; + \\
\; \\
\hspace{4em} \hat{u}_{x,1+1/2}v^{-}(x_{1+1/2})-\hat{u}_{x,1-1/2}v^{+}(x_{1-1/2})  \; +\\
\\
\hspace{2em}
\Bigg(  \dfrac{C_{1} }{h}     \left ( (u)^{+}_{1+1/2}-(u)^{-}_{1+1/2}   \right ) +
			 {C}_{2} \left (  (u_{x})^{+}_{1+1/2}  - (u_{x})^{-}_{1+1/2}  \right )
	  \Bigg)v^{-}_{1+1/2} \; -\\
\\
\hspace{2em}
\Bigg(  \dfrac{D_{1} }{h}  (u)^{+}_{1-1/2}+{D}_{2}(u_{x})^{+}_{1-1/2}\Bigg)v^{+}_{1-1/2} \; +\\
\\
\hspace{2em}
\Bigg( E_1    \left ( (u)^{+}_{1+1/2}-(u)^{-}_{1+1/2}   \right ) +
			 h{E}_{2} \left (  (u_{x})^{+}_{1+1/2}  - (u_{x})^{-}_{1+1/2}  \right )		
	\Bigg)(v_{x})_{1+1/2}^{-} \; -\\
\\
\hspace{2em}
\Bigg( F_{1}(u)^{+}_{1-1/2}+h F_{2}(u_{x})^{+}_{1-1/2}\Bigg)(v_{x})_{1-1/2}^{+}\nonumber
  \end{array}
  \right.
\end{equation}
\noindent where the general form of the fluxes are
\begin{equation}
\label{3_1}
\left.
\begin{array}{ll}
\hat{u}_{x,1+1/2}=\alpha(u_{x})^{-}_{1+1/2}+\left( 1-\alpha\right) (u_{x})^{+}_{1+1/2} \\
\\
\hat{u}_{x,1-1/2}=\beta(u_{x})^{-}_{1-1/2}+\left( 1-\beta\right) (u_{x})^{+}_{1-1/2}\\
  \end{array}
  \right.
\end{equation}
and $\alpha=1$. $\beta=0$ since the cell $I_{0}=\left[ x_{-{1/2}},x_{{1/2}}\right] $ does not exist.\\
\noindent Replacing $v(x)$ by $ \varphi_{1+1/4}$ and then $\varphi_{1-1/4}$ and comparing with Eq.(\ref{Chap3_1}), we obtain a system of eight equations for eight unknowns.  The solution to these equations is:\\

%

%
\begin{equation}
\label{Chap_coeff_solution_10} 
\left.
\begin{array}{llllll}
C_1=\frac{7}{3}, &&  C_2=\frac{1}{2},  \\
\\
D_1=\frac{14}{3}, &&  D_2=0, \\
\\
E_1= -\frac{1}{18} (8 c+5), && E_2=- \frac{1}{18} (c+1),\\
\\
F_1= \frac{1}{9}(8c+5), && F_2=0.\\

 \end{array} 
  \right.
\end{equation}

\noindent As for the second cell, the scheme is unchanged, but the stability is impacted by the scheme on the left side of the cell boundary at $x=h$.\\

\noindent The boundary operator ${\Theta}_{3/2}$ related to  $x_{3/2}=h$ can be adapted in the following way.\\

%
\begin{equation}
\label{Chap3_2}
\left.
\begin{array}{ll}
{\Theta}_{3/2}=
\bigg((u_{x})^{-}_{3/2}\bigg)u^{-}(x_{3/2})\\
-\bigg(\beta(u_{x})^{-}_{j-1/2}+\left( 1-\beta\right) (u_{x})^{+}_{j-1/2}\bigg)u^{+}(x_{3/2})\\
\\
+\bigg( -  \dfrac{C_{1}  }{h}  (u)^{-}_{3/2}-{C}_{2}(u_{x})^{-}_{3/2}\bigg)u^{-}_{3/2}\\
\\
-\Bigg( \dfrac{D_{1}  }{h}  \bigg((u)^{+}_{3/2}-(u)^{-}_{3/2}\bigg)+{D}_{2}\bigg((u_{x})^{+}_{3/2}-(u_{x})^{-}_{3/2}+\bigg)\Bigg)u^{+}_{3/2}\\
\\
+\bigg(-E_{1}(u)^{-}_{3/2}- h E_{2}(u_{x})^{-}_{3/2} \bigg)(u_{x})_{3/2}^{-}\\
\\
-\Bigg( F_{1}\bigg((u)^{+}_{3/2}-(u)^{-}_{3/2}\bigg)+ h F_{2}\bigg((u_{x})^{+}_{3/2}-(u_{x})^{-}_{3/2}\bigg)\Bigg)(u_{x})_{3/2}^{+}\\
\\
-\dfrac{1}{2}\int_{I_{1}}(u_{x})^2_{1}dx-\dfrac{1}{2}\int_{I_{2}}(u_{x})^2_{2}dx\nonumber
  \end{array}
  \right.
\end{equation}

\noindent where 

%
\begin{equation}
\label{Chap3_theta_32}
\left.
\begin{array}{llllll}
C_1=\frac{7}{3}, &&  C_2=\frac{1}{2},  \\
\\
D_1=\frac{7}{3}, &&  D_2=-\frac{1}{2}+\beta, \\
\\
E_1= -\frac{1}{18} (8 c+5), && E_2= -\frac{1}{18} (c+1),\\
\\
F_1= \frac{1}{18}(8c+5), && F_2=-\frac{1}{18}(c+1).\\

 \end{array} 
  \right.
\end{equation}

\noindent A sufficient condition for the stability is that ${\Theta}_{3/2}$ is non-positive definite.


\noindent The operator $\Theta_{3/2}$ may be written in the following form : 

\begin{equation}
\label{Chap3_2_1}
\left.
\begin{array}{ll}
{\Theta}_{3/2}=\sum_{k,l=1}^{4}a_{k,l}x_{k}x_{l}\nonumber
  \end{array}
  \right.
\end{equation}

\noindent where

\begin{equation}
\label{Chap3_2_2}
\left.
\begin{array}{ll}
\mathbf{x}=\begin{pmatrix} u_{3/4} \\  
   u_{5/4}\\
     u_{7/4}\\
      u_{9/4}\end{pmatrix}\nonumber
  \end{array}
  \right.
\end{equation}

\noindent The matrix for the associated symmetric bilinear form is
\begin{equation} \label{Theta_expansion}
 {\Theta}_{3/2}= 
    \mathbf{x}^{T} 
 M
  \mathbf{x}\nonumber
\end{equation}

is \\

\[M= \frac{1}{36h}
 \left[
    \scalemath{1}{
    \begin{array}{cccccccc}
     3(8c-19) & 71 - 40 c & -\frac{1}{2}(7-8c) & \frac{1}{2}(1-8c) \\
     \\
    71-40 c & 56 c-169 & \frac{1}{2} (113-40 c) & \frac{1}{2} (40 c-23) \\
    \\
     \frac{1}{2} (8 c-7) & \frac{1}{2} (113-40 c) & 56 c-169 & 71-40 c \\
     \\
     \frac{1}{2} (1-8 c) & \frac{1}{2} (40 c-23) & 71-40 c & 3 (8 c-19) \\
     \\
    \end{array}
    }
  \right]\nonumber
\]

\noindent It can be verified that for all values of $c$ between $-1$ and $1$, the eigenvalues of $36h*M$ are non-positive.

\subsubsection{Numerical Results}

\begin{figure}[h]
\begin{center}
\begin{tabular}{ccccc}
\hspace{-3.em}
\includegraphics[width=1.15\textwidth]{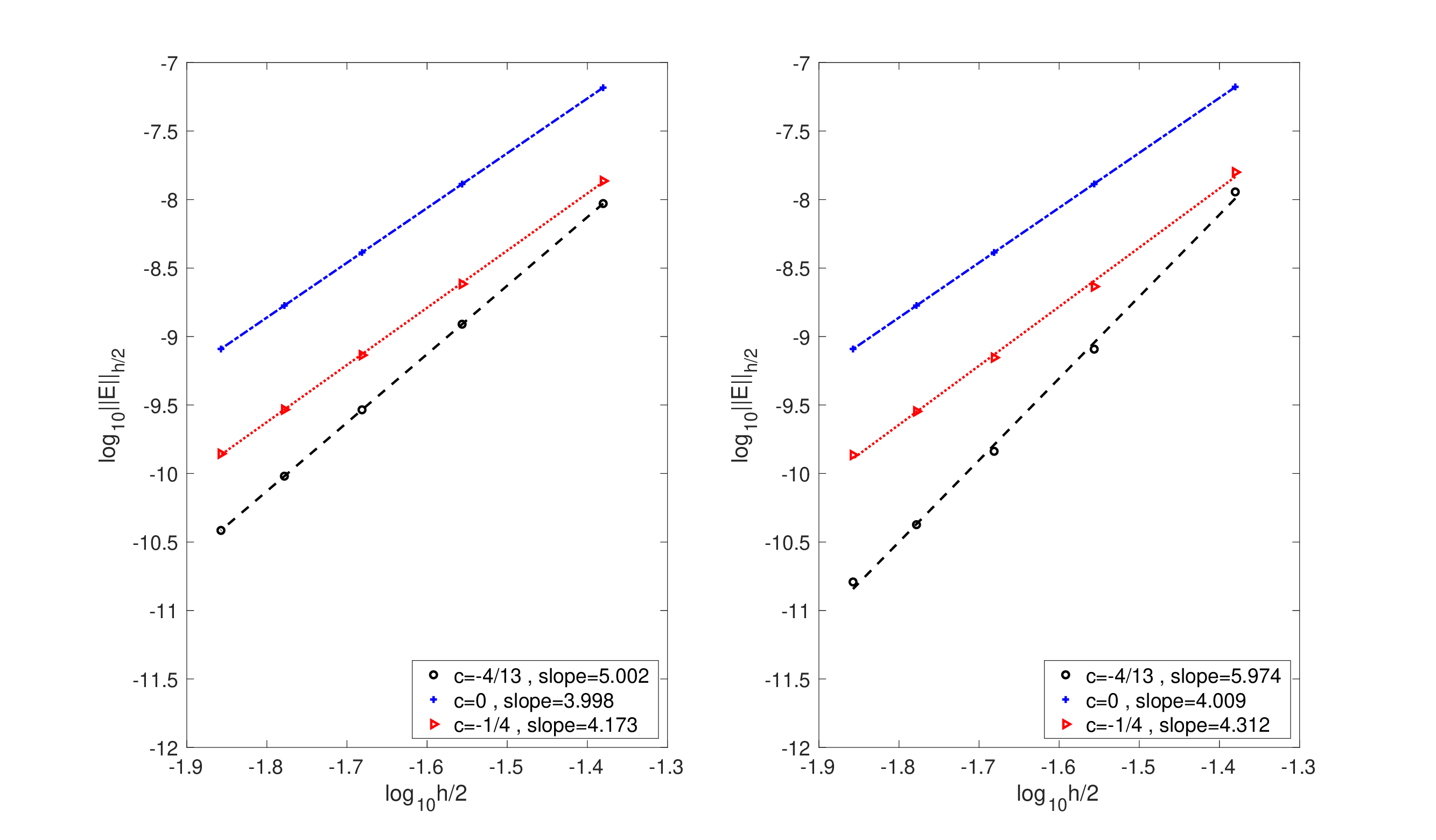}&&

\end{tabular}
\end{center}
\caption{Error and Convergence plots, $\log_{10} \|\vE\|$ vs. $\log_{10} h$ for 1D Heat Problem - Dirichlet boundary conditions. Left: without Post-Processing; right: with Post-Processing.
}\label{fig:fig_2}
\end{figure}

Numerical results demonstrate our theory on stability and order of convergence (see Fig. \ref{fig:fig_2}).
\noindent We used the approximation (\ref{dir3}) to solve the heat problem \eqref{21_dir} on the interval $[0,1]$ where the nonhomogeneous term, $F(x,t)$, and the initial condition were chosen such that the exact solution is $u(x,t)=\exp(\cos(x-t))$. The scheme was run with $N=24,36,48,60,72$. Sixth Order Explicit Runge-Kutta was used for time integration. This example shows that although this scheme has a third-order truncation error, this is a fourth-order scheme and fifth-order for the case of $c=-4/13$. \ref{app:Convergence} below shows that, in this case, an extra order can be achieved by using a post-processing filter. The right graph in Fig. \ref{fig:fig_2} verifies this analysis.

\subsection{Conclusion}\label{Conclusion}

We combined the advantages of viewing our BFD scheme first as a FE type and secondly as a FD type. This allowed us to provide a proof of stability in an easy manner based on the standard tools developed for DG methods, and then the optimal rate of convergence was proved through the tools developed for BFD methods.

\subsection{Generalization to Two Dimensions}

The generalization to two dimensions is fairly straightforward, as the two-dimensional scheme is constructed as a tensor product of the one-dimensional scheme presented in the previous subsection.


\noindent Let us first consider the following problem with periodic boundary conditions:
\begin{equation}
\label{21}
\left\lbrace 
\begin{array}{ll}
\dfrac{\partial u}{\partial t}(x,y,t)= \dfrac{\partial ^2}{\partial x^2 } u(x,y,t)+ \dfrac{\partial ^2}{\partial y^2 } u(x,y,t)+F(x,y,t)\mbox{ , } \\
\hspace{15 em}   (x,y)\in \left( 0,2\pi \right)\times \left( 0,2\pi \right)\mbox{ , }t\geq 0\\
u(x,y,t=0)=f(x,y)
\end{array}
  \right.
\end{equation}

\subsubsection{Proof of Optimal Convergence}

\begin{figure}[h]
\begin{center}
\begin{tabular}{ccccc}
\includegraphics[width=0.8\textwidth]{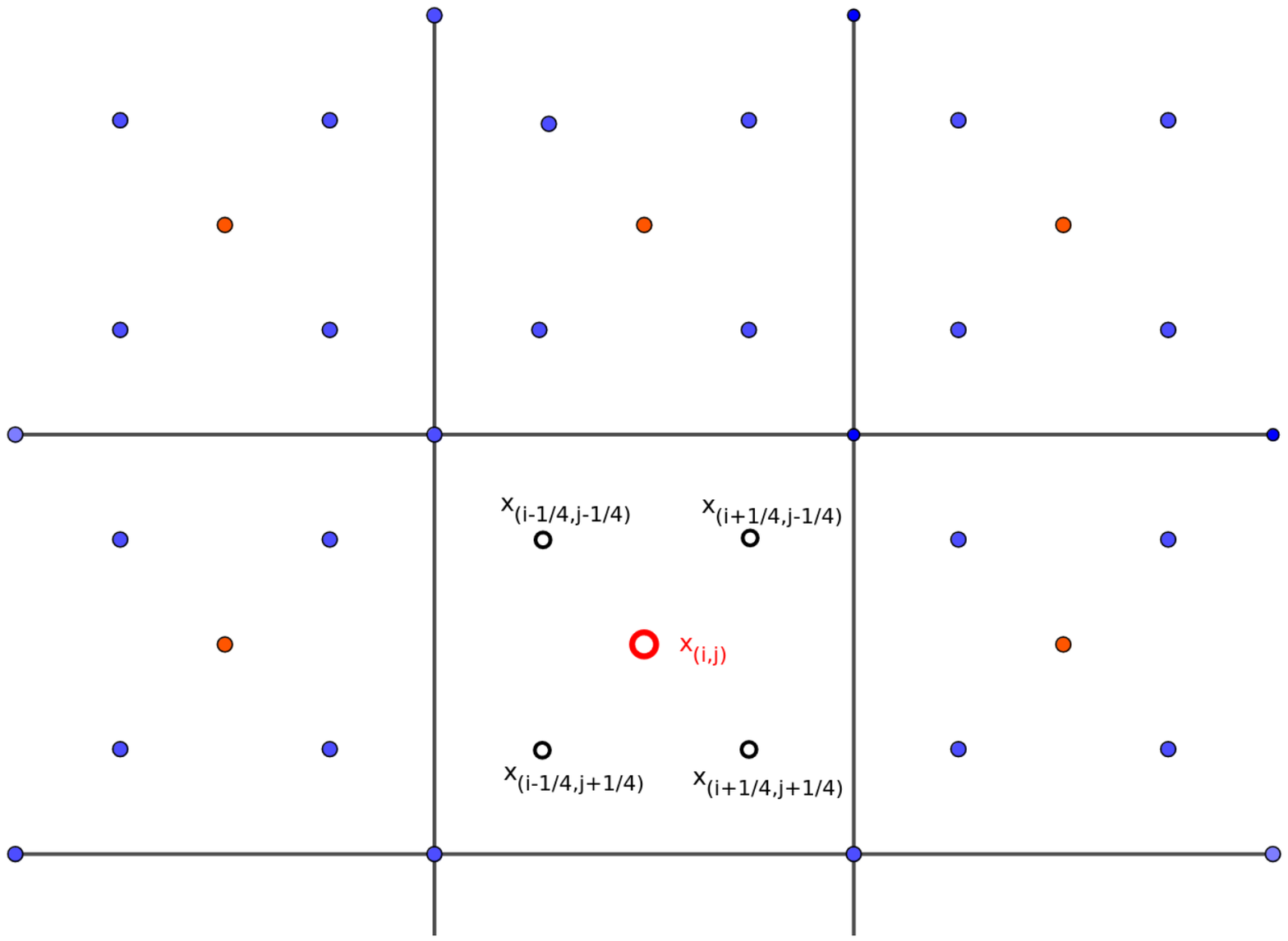}
\end{tabular}
\end{center}
\caption{Grid in 2 Dimensions - Illustration
}\label{fig:fig_101}
\end{figure}

Let us present the two-dimensional 5th order approximation of four-point block. We consider the following grid: In each cell whose center is located in $x_{(i,j)}$, we define four nodes as illustrated in Fig. \ref{fig:fig_101}:
\begin{eqnarray} \label{2D_grid}\nonumber 
\vx_{(i-{1/4},j-{1/4})} &=& \vx_{(i,j)} + \left(-\frac{h}{4},-\frac{h}{4}\right) \nonumber \\
\vx_{(i+1/4,j+1/4)} &=& \vx_{(i,j)} + \left(\frac{h}{4},\frac{h}{4}\right) \nonumber \\
\vx_{(i+1/4,j-1/4)} &=& \vx_{(i,j)} + \left(\frac{h}{4},-\frac{h}{4}\right) \nonumber \\
\vx_{(i-1/4,j+1/4)} &=& \vx_{(i,j)} + \left(-\frac{h}{4},\frac{h}{4}\right) \nonumber \\ 
&& \qquad h=2\pi/N , \; i,j=1,...,N  
\end{eqnarray}
\noindent where 
\begin{equation}\nonumber
\vx_{(i,j)}=\left(h(i-1)+\frac{h}{2},h(j-1)+\frac{h}{2}\right)
\end{equation}
\noindent Altogether there are $4N^2$ points on the grid, with a distance of $h/2$ between them when located at the same $x$ or $y$ coordinate.

\subsubsection{Proof of Stability}

The generalization to higher dimensions is done seamlessly. Indeed, the operator ${\Theta}_{j-1/2}$, as defined in the previous section for one dimension only, can be defined for two or more dimensions. In two dimensions, it is a linear combination of the contributions from each side of the cell and the average of two integrals of one dimension, in the $x$ and $y$ direction, respectively.\\

Similarly to the one-dimensional case, we may define all $\Theta$ operators in the x and y-directions respectively.\\

\begin{figure}[h]
\begin{center}	
\begin{tabular}{ccccc}
\includegraphics[width=1.5\textwidth]{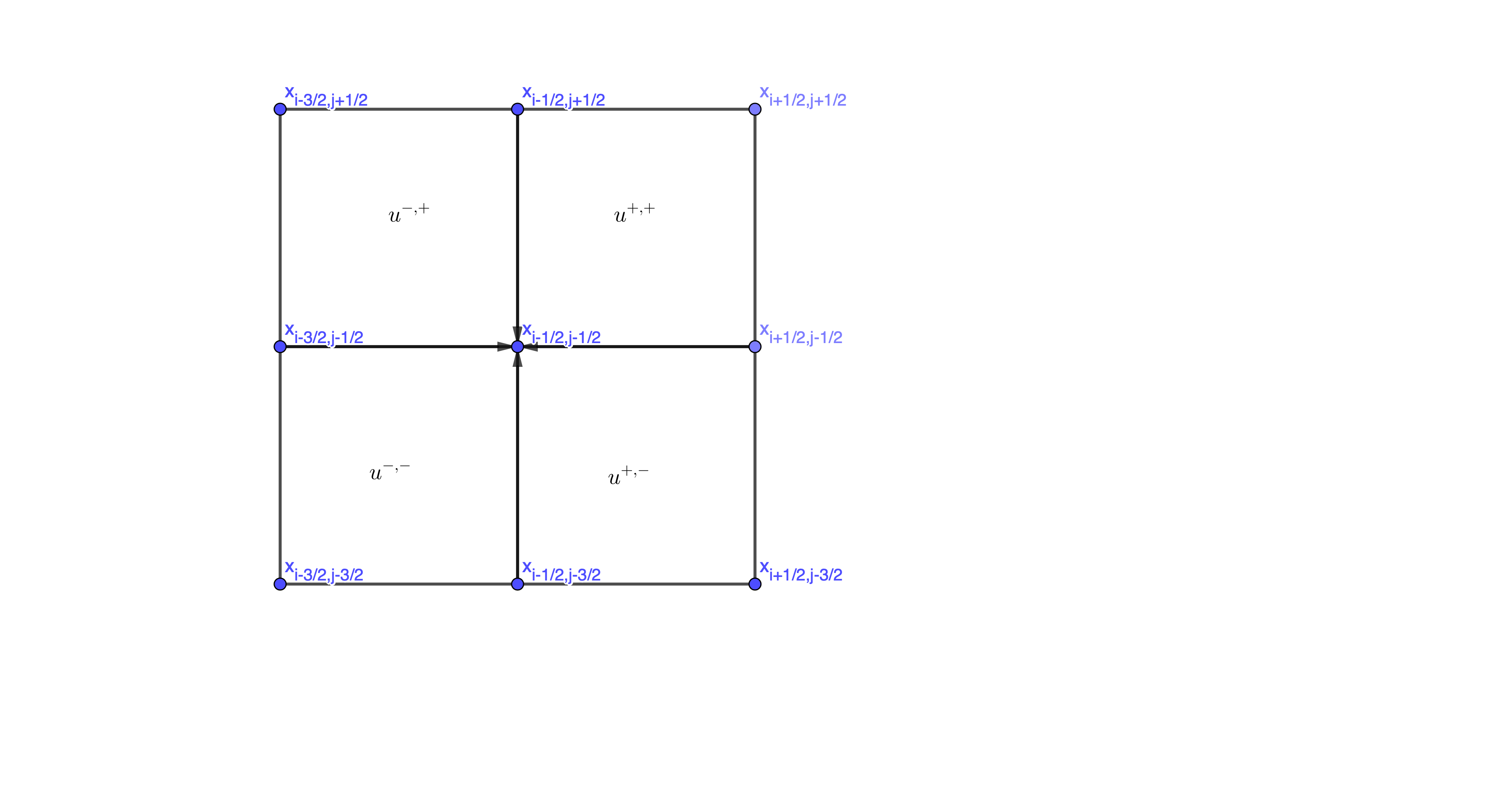}
\end{tabular}
\end{center}
\caption{Illustration of the dependence on node $\vx_{i-1/2,j-1/2}$ from neighbouring nodes in 2D
}\label{fig:fig_44}
\end{figure}

In two dimensions,  $\Theta$ is an integral over the edge of the cell. In a rectangular cell, four of those exist (see Fig. \ref{fig:fig_44}). Since our approximation space is a tensor product of polynomials degree up to  $k$, these integrals can be calculated using a convex combination of values of $u$ at the cell boundaries. In our case, $k=1$, and the linear combination is the average.\\

Hence, stability in the one-dimensional case implies stability in higher space dimensions.

\subsubsection{Numerical Results}

\begin{figure}[h!]
\begin{center}	
\begin{tabular}{ccccc}
\includegraphics[width=1\textwidth]{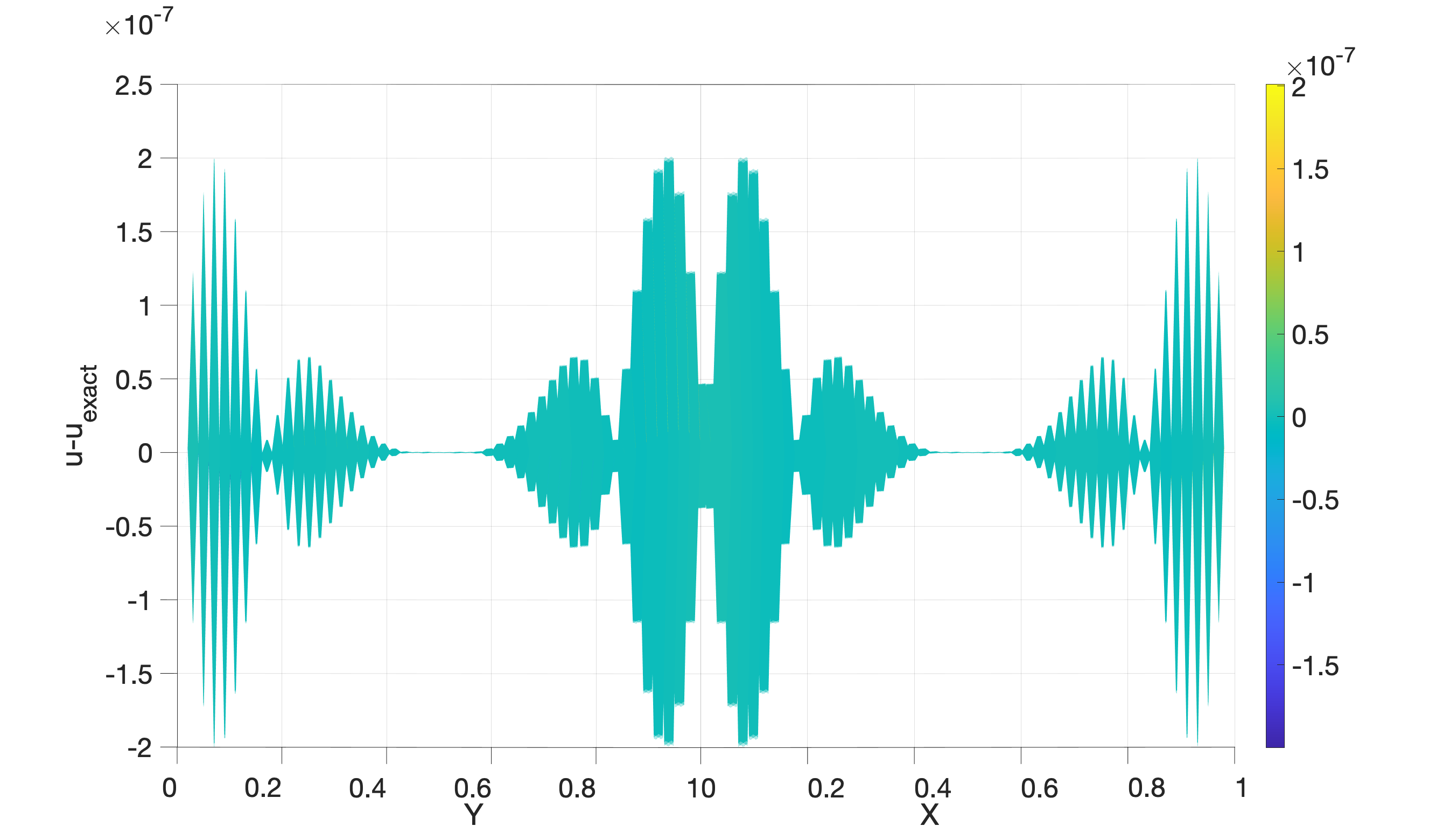}
\end{tabular}
\end{center}
\caption{2D Heat Problem - Two Points Block, BFD scheme - Periodic BC - Error at Final Time $T=1$ - $N=50$ - No post-processing
}\label{fig:fig_402}
\end{figure}

\begin{figure}[h!]
\begin{center}	
\begin{tabular}{ccccc}
\includegraphics[width=1\textwidth]{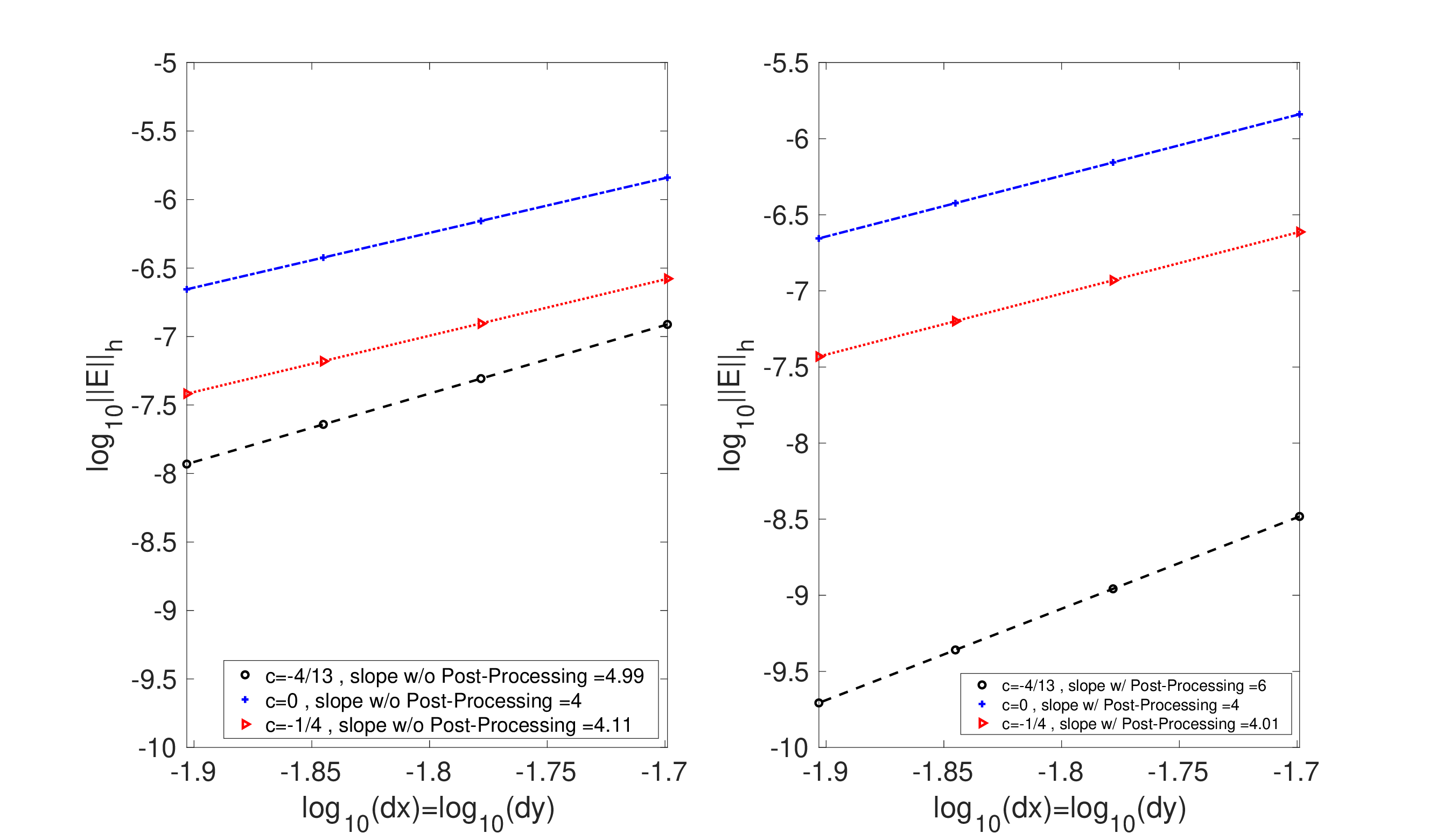}
\end{tabular}
\end{center}
\caption{2D Heat Problem - Two Points Block, BFD scheme, Convergence plots, $\log_{10} \|\vE\|$ vs. $\log_{10} h$ - Periodic BC - Left: no post-processing; right: Spectral post-processing
}\label{fig:fig_403}
\end{figure}

Numerical results demonstrate our theory on stability and order of convergence (see Fig. \ref{fig:fig_402}). The scheme was run for $u(x,y,t)=\exp(\cos(2\pi(x+y-t)))$ on the interval $[0,1]\times [0,1]$ and $N=50,60,70,80$ with a 6th order explicit Runge-Kutta time propagator and Final Time $T=1$. In order to clarify the error properties in 2D, the figure pictures the observation at $z=0 $ and far from $x=y=0$. In Fig. \ref{fig:fig_403}, the convergence plots appear with $c=-4/13,0,-1/4$, and the optimal convergence rate,  a fifth order without a  with post-processing and a sixth order with one is reached for $c=-4/13$.

\subsubsection{Non-Periodic Boundary Conditions}

Let us first consider the following two-dimensional Heat Problem with non-periodic boundary conditions:
\begin{equation}
\label{eq1_1}
\left\lbrace 
\begin{array}{ll}
\dfrac{\partial u}{\partial t}(x,y,t)= \dfrac{\partial ^2}{\partial x^2 } u(x,y,t)+ \dfrac{\partial ^2}{\partial y^2 } u(x,y,t)+F(x,y,t)\mbox{ , } \\
\hspace{15 em}   (x,y)\in \left( 0,\pi \right)\times \left( 0,\pi \right)\mbox{ , }t\geq 0\\
u(t=0)=f(x,y)\\
u(0,y,t)=g_{0}(y,t)\\
u(\pi,y,t)=g_{\pi}(y,t)\\
u(x,0,t)=h_{0}(x,t)\\
u(x,\pi,t)=h_{\pi}(x,t)
\end{array}
  \right.
\end{equation}

\begin{figure}[h]
\begin{center}

\includegraphics[width=1\textwidth]{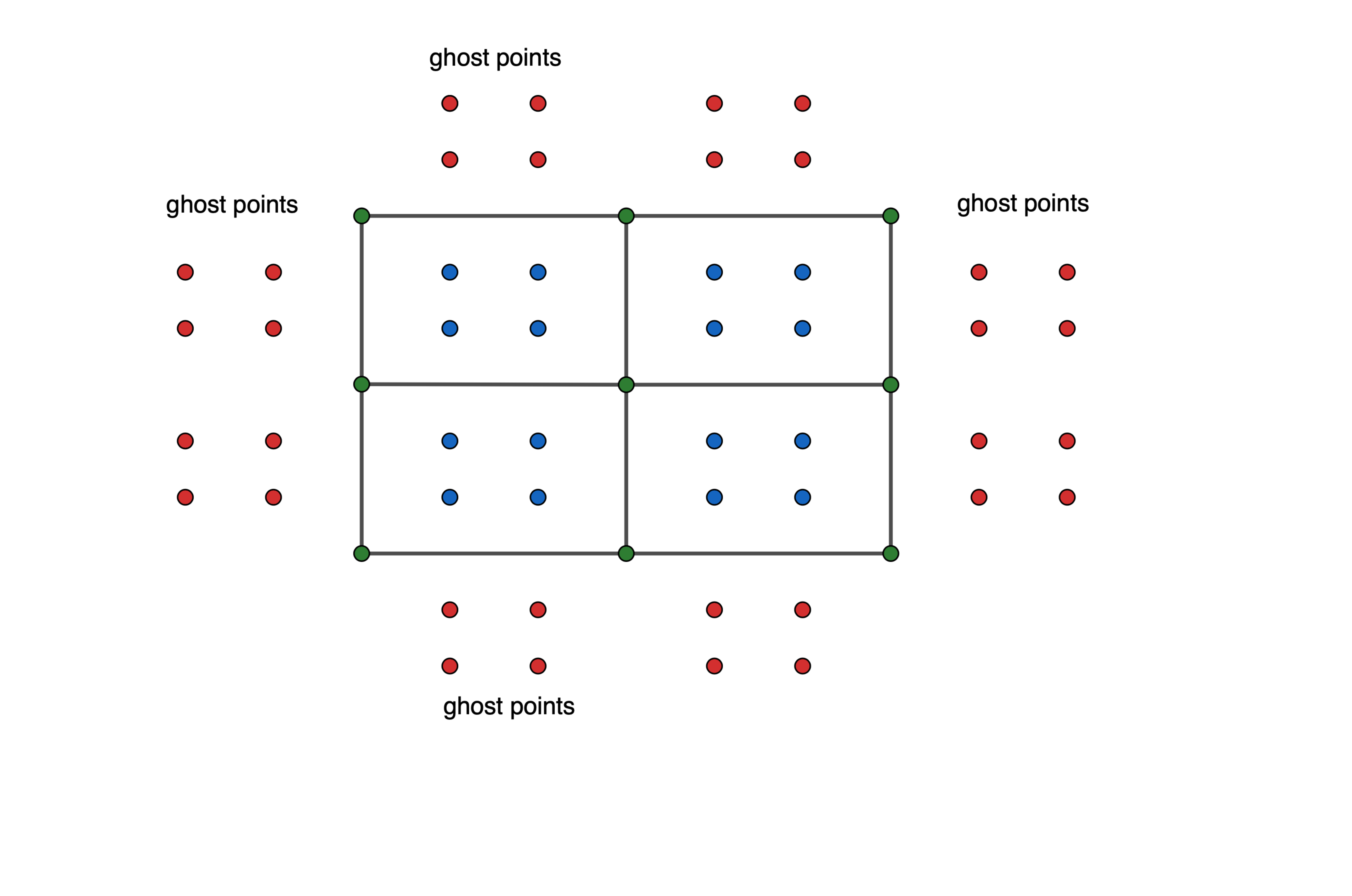}
\end{center}
\caption{Grid and Ghost Points in 2D - Illustration
} \label{fig:fig_5}
\end{figure}

\noindent In order to get an approximation scheme for the boundary points of the rectangle $[0,\pi]\times [0,\pi]$, we perform an extrapolation for the additional ghost points needed as illustrated in Fig.  \ref{fig:fig_5}.\\

\noindent In the x-axis direction:
$\vx_{(\frac{1}{4},j)}=(-\frac{h}{4},y_j)$, $\vx_{(-\frac{1}{4},j)}=(-\frac{3h}{4},y_{j})$ on one end\\
\\
and 
$\vx_{((N+1)-\frac{1}{4},j)}=(\pi+\frac{h}{4},y_{j})$, $\vx_{((N+1)+\frac{1}{4},j)}=(\pi+\frac{3h}{4},y_{j})$ 
on the other end.\\
\\
\noindent As for the internal points $\vx_{i,j}$, $1\leq i,j \leq N$, the scheme remains the same.

\noindent It follows that:

\begin{equation}
\label{dir11}
\left.
\begin{array}{ll}

u_{(-\frac{1}{4},j)}=-u_{(1+\frac{1}{4},j)}+2g_{0}(y_{j},t)+u_{xx}\left( 0,y_{j},t\right)\left( \frac{3h}{4}\right)^2+\\
\hspace{5em}  \frac{1}{12}u_{xxxx}\left( 0,y_{j},t\right)\left( \frac{3h}{4}\right)^4+O(h^6)\\

u_{(\frac{1}{4},j)}=-u_{(1-\frac{1}{4},j)}+2g_{0}(y_{j},t)+u_{xx}\left( 0,y_{j},t\right)\left( \frac{h}{4}\right)^2+\\
\hspace{5em}     \frac{1}{12}u_{xxxx}\left( 0,y_{j},t\right)\left( \frac{h}{4}\right)^4+O(h^6)\\

\end{array}
  \right.
\end{equation}

\noindent Similarly:

\begin{equation}
\label{dir22}
\left.
\begin{array}{ll}

u_{((N+1)-\frac{1}{4},j)}=-u_{(N+\frac{1}{4},j)}+2g_{\pi}(y_{j},t)+u_{xx}\left( \pi,y_{j},t\right)\left( \frac{h}{4}\right)^2+\\
\hspace{5em}   \frac{1}{12}u_{xxxx}\left( \pi,y_{j},t\right)\left( \frac{h}{4}\right)^4+O(h^6)\\

u_{((N+1)+\frac{1}{4},j)}=-u_{(N-\frac{1}{4},j)}+2g_{\pi}(y_{j},t)+u_{xx}\left( \pi,y_{j},t\right)\left( \frac{3h}{4}\right)^2+\\
\hspace{5em}   \frac{1}{12}u_{xxxx}\left( \pi,y_{j},t\right)\left( \frac{3h}{4}\right)^4+O(h^6)\\

\end{array}
  \right.
\end{equation}

\noindent where $u_{xx}(0,y_{i},t)$,$u_{xxxx}(0,t)$,$u_{xx}(\pi,y_{i},t)$,$u_{xxxx}(\pi,y_{i},t)$ can be directly computed from the PDE:

\begin{equation}
\label{}
\left.
\begin{array}{ll}

u_{xx}(0,y,t)=u_{t}(0,y,t)-u_{yy}(0,y,t)-F(0,y,t)\\
u_{xx}(\pi,y,t)=u_{t}(\pi,y,t)-u_{yy}(\pi,y,t)-F(\pi,y,t)\\
u_{xxxx}(0,y,t)=u_{tt}(0,y,t)-2u_{tyy}(0,y,t)-F_{t}(0,y,t)+\\
\hspace{8em}   u_{yyyy}(0,y,t)+F_{yy}(0,y,t)-F_{xx}(0,y,t)\\
u_{xxxx}(\pi,y,t)=u_{tt}(\pi,y,t)-2u_{tyy}(\pi,y,t)-F_{t}(\pi,y,t)+\\
\hspace{8em}   u_{yyyy}(\pi,y,t)+F_{yy}(\pi,y,t)-F_{xx}(0,y,t)\\

\end{array}
  \right.
\end{equation}

\noindent In the $y$-axis direction, the extrapolations are analogous.

\subsubsection{Numerical Results}

\begin{figure}[h!]
\begin{center}	
\begin{tabular}{ccccc}
\includegraphics[width=1\textwidth]{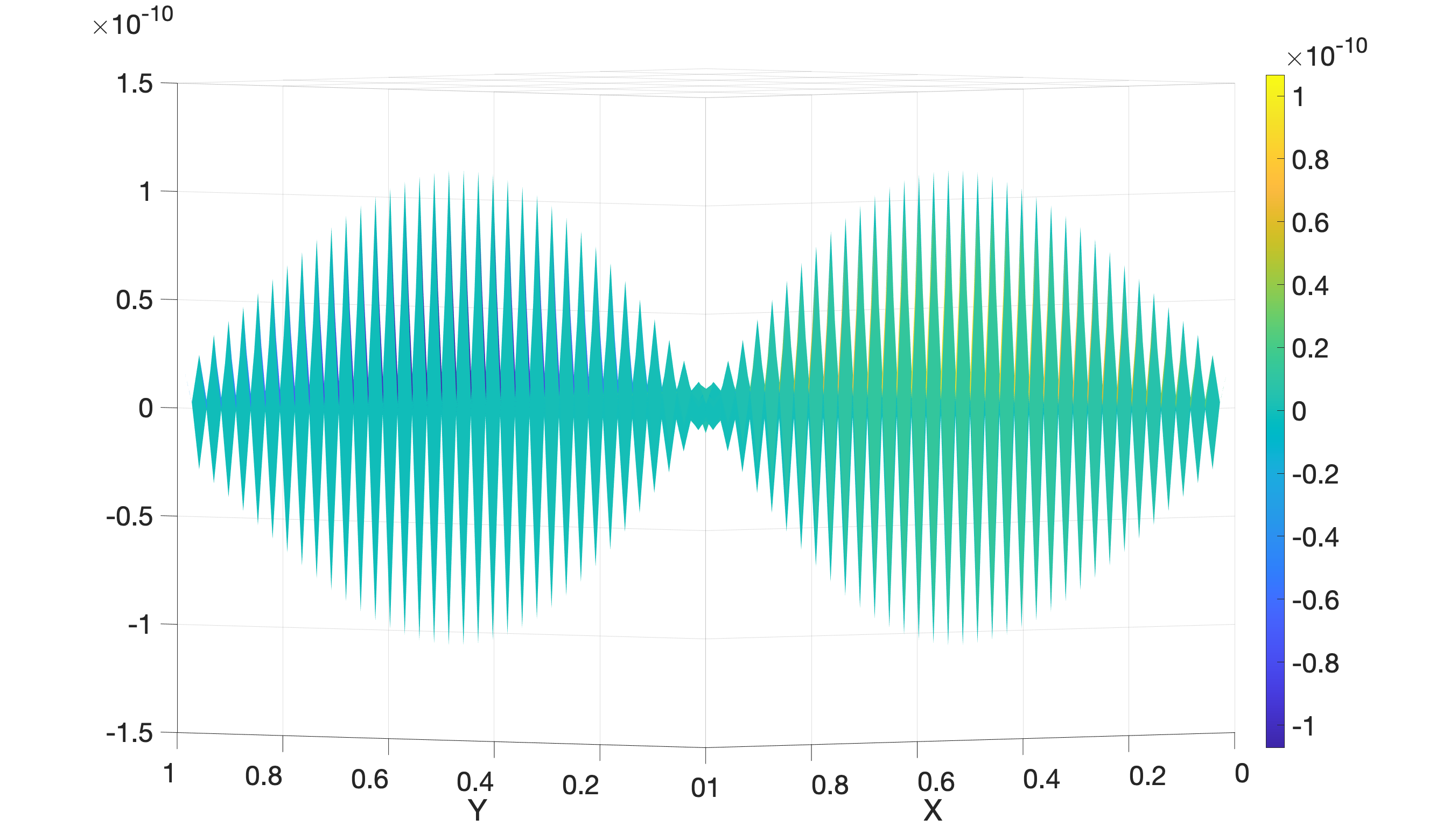}
\end{tabular}
\end{center}
\caption{2D Heat Problem - Two Points Block, BFD scheme - Dirichlet BC - Error at Final Time $T=1$ - $N=72$ - No post-processing
}\label{fig:fig_6}
\end{figure}

\begin{figure}[h!]
\begin{center}	
\begin{tabular}{ccccc}
\includegraphics[width=1\textwidth]{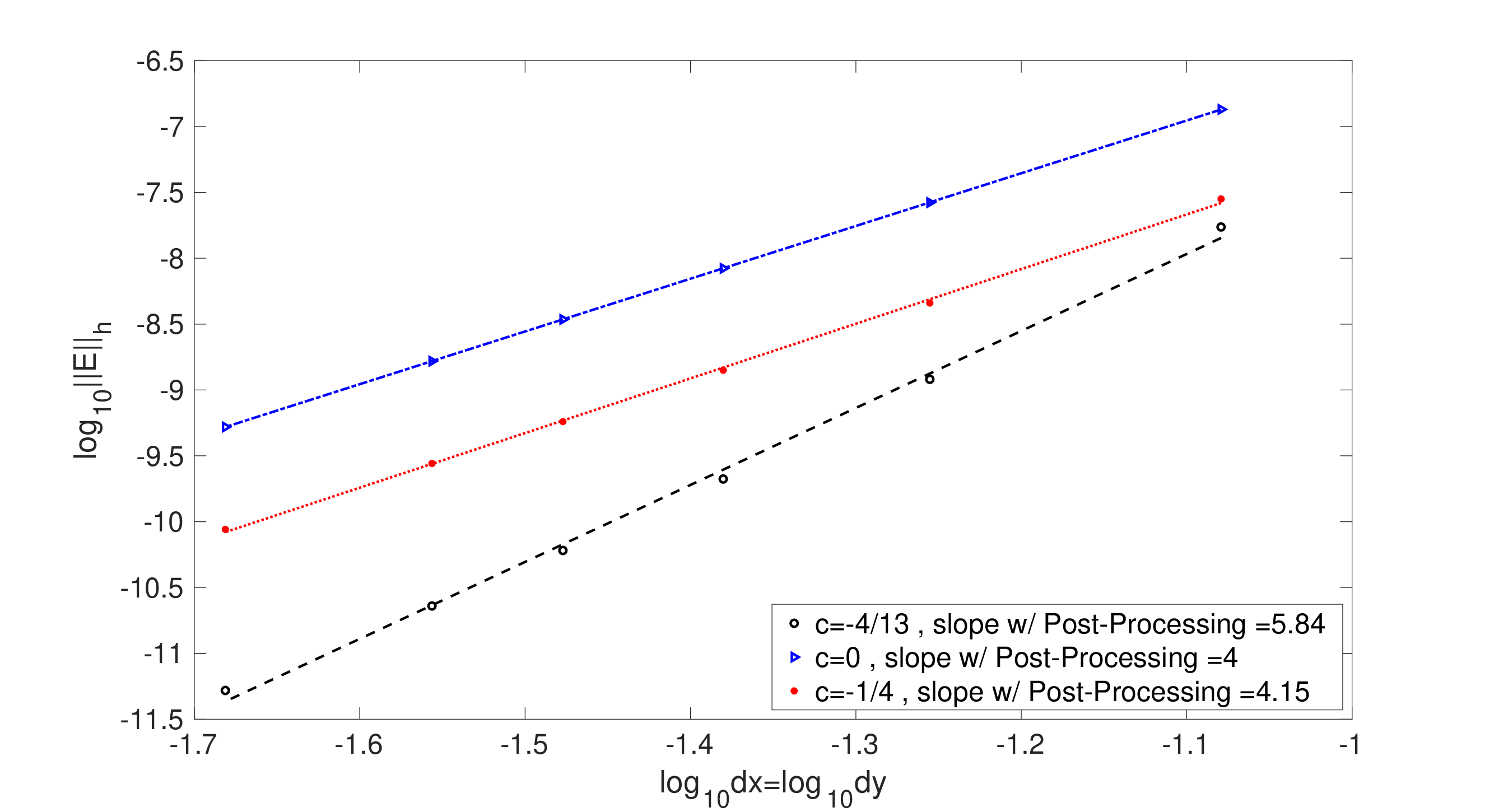}
\end{tabular}
\end{center}
\caption{2D Heat Problem - Two Points Block, BFD scheme, Convergence plots, $\log_{10} \|\vE\|$ vs. $\log_{10} h$ - Dirichlet BC -post-processing  filter
}\label{fig:fig_7}
\end{figure}

\noindent Numerical results demonstrate our theory on stability and order of convergence (see Fig. \ref{fig:fig_6}). The scheme was run for an exact solution $u(x,y,t)=\exp(\cos(x+y-t))$ on the interval $[0,1]\times [0,1]$ and $N=24,36,48,60,72,96$ with a sixth order explicit Runge-Kutta time propagator and Final Time $T=1$. In order to clarify the error properties in 2D, the figure pictures the view from an observer at $z=0$ and far from $x=y=0$.\\

\noindent  The convergence plots present for $c=0,-1/4,-4/13$ BFD schemes, are presented in Fig. \ref{fig:fig_7}. The same optimal convergence rate is reached for $c=-4/13$, resulting in a sixth-order convergence rate with a post-processing filter. The scheme was run for an exact solution $u(x,y,t)=\exp(\cos(x+y-t))$ on the interval $[0,1]\times [0,1]$ and $N=24,36,48,60,72,96$ with a sixth-order explicit Runge-Kutta time propagator and Final Time $T=1$.

\subsection{Generalization to Three Dimensions}

In this section, we present a brief description of the three-dimensional BFS/GD scheme.


\noindent Consider the Heat equation  with periodic boundary conditions:
\begin{equation}
\label{3D_10}
\left\lbrace 
\begin{array}{ll}
\dfrac{\partial u}{\partial t}(x,y,zt)=\nabla^2 u(x,y,z,t) )+F(x,y,t)\mbox{ , } \\
\hspace{10 em} \vx=  (x,y,z)\in \left( 0,2\pi \right)\times \left( 0,2\pi \right)\times \left( 0,2\pi \right)\mbox{ , }t\geq 0\\
u(x,y,z,t=0)=f(x,y,z)
\end{array}
  \right.
\end{equation}  
The domain is divided into cells, $I_{(j,k,l)}$, where 
\begin{equation}  \label{3D_20}
I_{(j,k,l)} = \left \{ \vx : \left | x-x_{(j,k,l)} \right | ,   \left | y-y_{(j,k,l)} \right | ,  \left | z-z_{(j,k,l)} \right |  \le \frac{h}{2}\right \}
\end{equation}
and
\begin{equation}    \label{3D_30}
\vx_{(j,k,l)}= \left ( x_{(j,k,l)}, y_{(j,k,l)}, z_{(j,k,l)}\right ) = \left(h(j-1)+\frac{h}{2}, h(k-1)+\frac{h}{2} ,  h(l-1)+\frac{h}{2}  \right)
\end{equation}
The grid points, at each cell, are
\begin{equation}    \label{3D_40}
\vx_{(j \pm 1/4 ,k \pm 1/4 ,l \pm 1/4 )} = \vx_{(j,k,l)} + \left ( \pm 1, \pm 1, \pm 1 \right )\frac{h}{4} 
\end{equation}
Altogether, there are eight grid points at each cell. 

These nodes form a rectangular grid aligned with the axes. The one-dimensional BFD scheme, \eqref{2.6}, is used for each line parallel to the $x$, $y$, and $z$. This three-dimensional scheme is stable since the operator $\Theta$ at each face of the cells $I_{(j,k,l)}$ is a convex combination of the non-positive one-dimensional $\Theta$. Thus, it is non-positive as well.\\

 For the non-periodic boundary conditions, as for the two-dimensional case, to approximate the values at the boundary points of the cube $[0,\pi]\times [0,\pi]\times [0,\pi]$, we perform an extrapolation to obtain the additional ghost points. The formulation for these extrapolations and the proofs are similar to the one- and two-dimensional cases. \\

\section{Future Work}\label{Future Work}

The schemes we constructed and analyzed in this manuscript are defined in rectangular domains. Our project's next stage is to derive block finite difference schemes to solve the heat equation in two and three-dimensional, complicated geometries.

We also intend to derive block finite difference schemes to solve the advection equation and hyperbolic systems.

\section{Acknowledgement}
This research was supported by Binational (US-Israel) Science Foundation grant No. 2016197.



 \bibliographystyle{elsarticle-num} 
\bibliography{ref_tmp}

\appendix

\section{Pseudo-Fourier Analysis - Eigenvalues and Eigenvectors} \label{app:Convergence}

\noindent This section will use tools developed in \citep{ditkowski2015high} and \citep{ditkowski2020error}. The goal of the procedure is to find a value of the parameter $c$ that appears in the BFD scheme \eqref{BFD_scheme_10}, leading to the optimal convergence rate of this scheme.
\\
\noindent We will perform a Fourier-like analysis for this purpose. In \citep{GUO2013458}, a straightforward Fourier analysis was performed using different equations with periodic boundary conditions. Here, the analysis of the error structure basis is a linear combination of eigenvectors, hence covering the case where the matrix $Q$ is not diagonalizable using standard discrete Fourier transform. The analysis still requires periodic boundary conditions and a uniform mesh.\\

\noindent We first split the Fourier spectrum into low and high frequencies as follows:
Let $\omega \in \{-N/2,\ldots,N/2$\} being an integer and :
\begin{equation}\label{3.30} 
\nu \,=\,   \left\{ \begin{array}{lcl}
                       \omega -N,\hspace{1cm}&\omega >0\\
                       \omega+N,\hspace{1cm} &\omega\leq 0
                     \end{array}
\right.
\end{equation}
\noindent Then, 
\begin{eqnarray}   \label{3.40} 
&\text{for } \omega \ge0:& e^{i \nu x_{j-1/4}} = -i e ^{i  \omega x_{j-1/4}}  \;\; {\rm and } \;\; 
 			e^{i \nu x_{j+1/4}} = i e ^{i  \omega x_{j+1/4}}  \nonumber \\
&\text{for } \omega < 0:& e^{i \nu x_{j-1/4}} = i e ^{i  \omega x_{j-1/4}}  \;\; {\rm and } \;\; 
 			e^{i \nu x_{j+1/4}} = -i e ^{i  \omega x_{j+1/4}} 
\end{eqnarray}

\noindent We now present thae analysis for $\omega\geq0$. The analysis for  $\omega < 0$ is similar.


\noindent We denote the vectors ${\rm e}^{i \omega \bm{x}}$ and ${\rm e}^{i \nu \bm{x}}$ by:

\begin{equation}\label{3.50} 
{\rm e}^{i \omega \bm{x}} =  \left(
                          \begin{array}{c}
                            \vdots \\
                            {\rm e}^{i \omega x_{j-{1/4}}} \\
                            {\rm e}^{i \omega x_{j+1/4}} \\
                            \vdots
                          \end{array}
                        \right) \mbox{, }{\rm e}^{i \nu \bm{x}}= \left(
                          \begin{array}{c}
                            \vdots \\
                            {\rm e}^{i \nu x_{j-{1/4}}} \\
                            {\rm e}^{i \nu x_{j+1/4}} \\
                            \vdots
                          \end{array}
                        \right)
\end{equation}

\noindent We look for eigenvectors in the form of:

\begin{equation}\label{3.50} 
\psi_k(\omega) ={\alpha_k}   \frac{   {\rm e}^{i \omega \bm{x}}  }{\sqrt{2 \pi}}    +
{\beta_k}    \frac{  {\rm e}^{i \nu \bm{x}} }{\sqrt{2 \pi}} 
  \end{equation}
                      
\noindent where, for normalization, it is required that
$|\alpha_k|^2+|\beta_k|^2=1 $, $k=1,2$. 

\noindent We note here that each component of the linear combination ${\rm e}^{i \omega \bm{x}}$ and ${\rm e}^{i  \nu \bm{x}}$ is not an eigenvector in the classical sense, since only the linear combination verifies the initial condition.

\noindent However:

\begin{equation}\label{3.50} 
Q{\rm e}^{i \omega \bm{x}} ={\rm diag}({\mu}_{1},{\mu}_{2},...,{\mu}_{1},{\mu}_{2}){\rm e}^{i \omega \bm{x}}
\end{equation}

\begin{equation}\label{3.501} 
Q{\rm e}^{i \nu\bm{x}} ={\rm diag}({\sigma}_{1},{\sigma}_{2},...,{\sigma}_{1},{\sigma}_{2}){\rm e}^{i \nu \bm{x}}
\end{equation}

\noindent where:
\begin{eqnarray}  \label{3.51_mu_sigma} 
\mu_1  &=& \frac{8}{3h^2}\left [
					-  { \sin ^2 \left( \frac{h  \omega}{4}\right)  
								\left(7-\cos \left(\frac{h  \omega}{2} \right)  \right)}
					- 4 \, i\,  c \,  e^{\frac{i h w}{4}} \sin ^5\left(\frac{h w}{4}\right) 
				 \right ]  \nonumber \\ 
				 \nonumber \\
\mu_2  &=&   \frac{8}{3h^2}\left [
					-  { \sin ^2 \left( \frac{h  \omega}{4}\right)  
								\left(7-\cos \left(\frac{h  \omega}{2} \right)  \right)}
					+  4 \, i\,  c \,  e^{-\frac{i h w}{4}} \sin ^5\left(\frac{h w}{4}\right) 
				 \right ]  \nonumber \\ 
				 \nonumber \\
\sigma_1  &=& \frac{8}{3h^2}\left [
					{ -\cos ^2\left(\frac{h  \omega}{4}\right) 
								\left(7+\cos \left(\frac{h  \omega}{2} \right)  \right)}
					+  4 \ c \,  e^{\frac{i h w}{4}} \cos^5\left(\frac{h w}{4}\right) 
				 \right ]  \nonumber \\ 
				 \nonumber \\
\sigma_2  &=&  \frac{8}{3h^2}\left [
					{ -\cos ^2\left(\frac{h  \omega}{4}\right) 
								\left(7+\cos \left(\frac{h  \omega}{2} \right)  \right)}
					+  4 \ c \,  e^{-\frac{i h w}{4}} \cos^5\left(\frac{h w}{4}\right) 
				 \right ]  \nonumber 
				 \nonumber \\
\end{eqnarray}

\noindent In order to find the coefficients ${\alpha}_{k}$ and ${\beta}_{k}$ along with the eigenvalues (symbols) $\hat{Q}_k$ for $\omega > 0$, we consider some node $x_{j}$.
\noindent Then we solve the following system of equations:

\begin{equation}
\label{3.512}
\left. 
\begin{array}{ll}
{\mu}_{1}\frac{\alpha_k}{\sqrt{2 \pi}} {\rm e}^{i \omega {x_{j-\frac{1}{4}}}} +{\sigma}_{1}\frac{\beta_k}{\sqrt{2 \pi}} {\rm e}^{i \nu {x_{j-\frac{1}{4}}}}=
{\hat{Q}}_k\left( \frac{\alpha_k}{\sqrt{2 \pi}} {\rm e}^{i \omega {x_{j-\frac{1}{4}}}} +\frac{\beta_k}{\sqrt{2 \pi}} {\rm e}^{i \nu {x_{j-\frac{1}{4}}}}\right)\\
{\mu}_{2}\frac{\alpha_k}{\sqrt{2 \pi}} {\rm e}^{i \omega {x_{j+\frac{1}{4}}}} +{\sigma}_{2}\frac{\beta_k}{\sqrt{2 \pi}} {\rm e}^{i \nu {x_{j+\frac{1}{4}}}}=
{\hat{Q}}_k\left( \frac{\alpha_k}{\sqrt{2 \pi}} {\rm e}^{i \omega {x_{j+\frac{1}{4}}}} +\frac{\beta_k}{\sqrt{2 \pi}} {\rm e}^{i \nu {x_{j+\frac{1}{4}}}}\right)
\end{array}
  \right.
\end{equation}

\noindent We denote $r_{k}=i\dfrac{{\beta}_{k}}{{\alpha}_{k}}$ and use Eq. \eqref{3.40} to obtain a simpler system:

\begin{equation}
\label{3.513}
\left. 
\begin{array}{ll}
{\mu}_{1}-{\sigma}_{1}r_{k}=
{\hat{Q}}_k\left( 1-r_{k}\right)\\
{\mu}_{2}+{\sigma}_{2}r_{k}=
{\hat{Q}}_k\left( 1+r_{k}\right)\\
\end{array}
  \right.
\end{equation}
\noindent Consequently, the expressions for
$r_k$ and the eigenvalues (symbols) $\hat{Q}_k$ are:

\begin{equation}
\label{3.513}
\left. 
\begin{array}{ll}
r_{1}=\dfrac{ \Omega+\Delta}{16 c \sin\left(\frac{h\omega}{4}\right)\cos ^5\left( \frac{h\omega}{4}\right)}i \\\\
r_{2}=\dfrac{ \Omega-\Delta}{16 c \sin\left(\frac{h\omega }{4}\right)\cos ^5\left( \frac{h\omega}{4}\right)}i \\\\
\end{array}
  \right.
\end{equation}
\noindent where
\begin{eqnarray}\label{3.56} 
\Omega  &=&  -\cos\left( \frac{h\omega }{2}\right)
				\Big( 16- \Big( 7+\cos\left( h\omega \right) \Big ) c \Big ) \\
\text{and} \nonumber  \\
\Delta  &=& \sqrt{\Omega^2 + 4 c^2  \sin^6 \left( \frac{h\omega }{2}\right)}
\end{eqnarray}
The symbols $\hat{Q}_{1}(\omega)$ and $\hat{Q}_{2}(\omega) $ are:
\begin{equation}\label{3.58} 
\left. 
\begin{array}{ll}
\hat{Q}_{1}(\omega)  = \dfrac{2}{3 h^2} \Big (- \Big (15+ \cos (h \omega )  \Big )  +  \Big (5+3 \cos (h \omega )  \Big ) c  + \Delta\Big)     \mbox{, } \\
\\
\hat{Q}_{2}(\omega)  = \dfrac{2}{3 h^2} \Big (- \Big (15+ \cos (h \omega )  \Big )  +  \Big (5+3 \cos (h \omega )  \Big ) c  + \Delta\Big) \\
\end{array}
  \right.
\end{equation}
\noindent Using the normalization condition $|\alpha_k|^2+|\beta_k|^2=1 $, $k=1,2$, we choose the coefficients ${\alpha}_{k},{\beta}_{k}$ to be:

\begin{equation}\label{3.57} 
\alpha_{1}=\dfrac{1}{\sqrt{1+{\vert r_{1}\vert}^2}}\mbox{, }\beta_{1}=-i\dfrac{r_1}{\sqrt{1+{\vert r_{1}\vert}^2}}\\
\end{equation}

\begin{equation}\label{3.58} 
\alpha_{2}=i \dfrac{ r_2 /\vert r_{2}\vert }{\sqrt{1+{\vert r_{2}\vert}^2}}\mbox{, }\beta_{2}=\dfrac{\vert r_2\vert}{\sqrt{1+{\vert r_{2}\vert}^2}}
\end{equation}
It is important to note here that all eigenvalues $\hat{Q}_1$  and $\hat{Q}_2$  are real and non-positive for all $\omega$ and $|c| \le 1$. Therefore, the scheme is Von-Neumann
stable. Additionally, since a complete set of eigenvectors exists, we conclude that the scheme is stable.

\bigskip

\noindent For $\omega h \ll 1$ the eigenvalues are:
\begin{equation}\label{3.60} 
\hat{Q}_1(\omega) =-\omega^{2} 
		+\dfrac{(4+13c)\omega^6{ h}^4}{2880(2-c)} 
		-\dfrac{(4 + 38 c + c^2) \omega ^8  h^6}{64512 (2-c)^2}
		+O(h^8)\\
\end{equation}
\noindent and
\begin{equation}\label{3.80} 
\hat{Q}_2(\omega) =-\dfrac{32(2-c)}{3h^2}
		+\dfrac{(5-6c)\omega^2}{3}
		-\dfrac{(1-3c)\omega^4 h^2}{18}
		+O(h^4)
\end{equation}
Since $\hat{Q}_2(\omega) = O(1/h^2)$ is negative, it immediately decays. Therefore ,we can neglect the $\psi_2$ terms.
\bigskip

If the initial condition is $u(x, t=0)={\rm e}^{i \omega x}/\sqrt{2 \pi}$, i.e. 
\begin{equation}\label{3.100} 
\vv_{j-\frac{1}{4}}(0) = {\rm e}^{i \omega x_{j-\frac{1}{4}}}/\sqrt{2 \pi} \;,
		 \vv_{j+\frac{1}{4}}(0) = {\rm e}^{i \omega x_{j+\frac{1}{4}}}/\sqrt{2 \pi}
		  \;\; ; \; \;\;\; \omega^2h\ll 1
\end{equation}
\noindent Then

\begin{equation}\label{3.110} 
\left. 
\begin{array}{ll}
{\vv}_{j-\frac{1}{4}}(t) =
\left( {\rm e}^{-\omega^2 t}\left[ 1+
		\dfrac{(4+13c){\omega}^{6} {h}^{4}}{2880(2-c)} t \right]+O(h^6)\right) 
		 \dfrac{  {\rm e}^{i \omega x_{j-\frac{1}{4}}} }{\sqrt{2 \pi}} +\\
\\
\hspace{8em} \left(\dfrac{c \, {\rm e}^{-\omega^2 t}({\omega h})^{5}}{1024(2-c)}+O(h^6)\right)  				\dfrac{   {\rm e}^{i \nu x_{j-\frac{1}{4}}}  }{\sqrt{2 \pi}} 
\end{array}
  \right.
\end{equation}
The same expression holds for $x_{j+1/4}$. The scheme is of fourth-order in general. However, if $c=-4/13$, then it becomes fifth-order. Furthermore, upon closer inspection of equation \eqref{3.110}, it appears that the error term's leading coefficient is of high frequency for the $c=-4/13$ case. This suggests that applying a post-processing filter can remove this high-frequency error, resulting in a sixth-order convergence rate.

\section{Post-processing}\label{app:Post-processing}


The analysis presented in \ref{app:Convergence} demonstrates that when the initial condition is $u(x, t=0)= \exp( i \omega x)$ and $\omega h \ll 1$, the numerical solution at time $t_n$ can be calculated using equation \eqref{3.110}. This scheme is of fourth-order, but selecting $c=-4/13$ eliminates the leading term in error, resulting in a fifth-order convergence rate. We now observe that the fifth-order term is of high frequency; therefore, it can be filtered in a post-processing stage to get a sixth-order scheme. This section describes the construction of these post-processing filters.

\subsection{Post-Processing applied to Heat problem with Periodic Boundary Conditions}

The analysis in \ref{app:Convergence} showed that the sixth-order convergence rate could be recovered by choosing $c=-4/13$ and filtering the high-frequency Fourier modes (larger than $N/2$ in absolute value) at a post-processing stage.

Thus, for the one-dimensional problem, the post-processing filtering is done as follows:
\begin{enumerate}
\item Computing the discrete Fourier transform (DFT) of the numerical solution at the final time, $\vv^n$.
\item Zeroing the coefficients of the high-frequency modes.
\item  Computing the inverse discrete Fourier transform.
\end{enumerate}

In the two-dimensional case, this procedure is performed for every row and column.

\bigskip


\subsection{Post-Processing for Non-periodic boundary conditions of Dirichlet type}


We look for a one-step procedure leading to an increased order of accuracy.
\noindent Unlike the case of periodic boundary conditions, the filter presented above cannot be used. Hence, a different technique shall be applied. Those filters were derived from a class of filters specifically applied to DG methods in the case of hyperbolic problems (see \cite{10.1007/978-3-319-19800-2_6}, \cite{ryan2005}, \cite{cockburn2003enhanced}), and compared with the standard \textit{Savitzky-Golay Filter} \cite{savitzky}.

\subsubsection{Filter Description}

This section is dedicated to the description of the three filters that were used, starting with the standard filter \textit{Savitzky-Golay Filter} \cite{savitzky}.\\

\noindent \textbf{Savitzky-Golay Filter}\\

Since their introduction more than half a century ago, Savitzky–Golay (SG) filters have been popular in many fields of data processing. They use a method of data smoothing based on local least-squares polynomial approximation. SG filters are usually applied to equidistant data points and are based on fitting a polynomial of a given degree $n$ to the data in a (usually symmetric) neighborhood $k-m...k+m$ of each data point $k$. This range contains $2m + 1$ data points. Each data point is replaced by the value of the fit polynomial at this point $k$. As this process is a linear filter and takes a limited number of points as the input, SG smoothing belongs to the class of Finite Impulse Response (FIR) filters. Therefore, it can be implemented as a convolution with a suitable kernel.\\

The filter was used as a benchmark for our custom filter, with the parameters $n=6,m=5$, using Matlab\textcopyright \cite{MATLAB:2021}  \textit{sgolayfilt} tool.
\\
We note here that $m$ is limited to odd values in order to conserve the symmetry of the interpolation. Asymmetric filters can be designed, but their efficiency does not increase \cite{luo}.

Numerical results will show that this filter is efficient from the error reduction aspect, but does not affect the order of accuracy of the BFD scheme.

In the next paragraph, we try and design filters based on a different kernel.\\

\noindent \textbf{First Interpolation Filter}\\

\noindent For every set of the 12 consecutive grid points, we perform an interpolation procedure of the approximated solution valued at those grid nodes to obtain a sixth-order interpolation polynomial (for the optimal $c=-4/13$ case).\\

\noindent We then compute the related approximation for the set of points and proceed to the next set of 12 nodes. This is in contrast with the SG filter, where each node is approximated by a different polynomial.

\noindent Since the post-processed error is still somewhat oscillatory, a second version of an interpolation filter is designed in the following paragraph.\\

\noindent \textbf{Second Interpolation Filter}\\

\noindent This time, the first six points from the left are interpolated as before. Then, the filtered values of each set of two points are taken as the six-order interpolation based on these two points, five points from the left and five from the right. The last six points are computed similarly to the first six.
There is a significant improvement in the Final Time error, as illustrated in Fig. \ref{fig:fig_post_005}. A seventh order of convergence is reached.\\

\subsubsection{Post-Processing Numerical Results}


\noindent  We employed the approximation \eqref{dir3}-\eqref{dir5} to find the solution to the heat problem \eqref{21_dir} on the interval $[0,1]$. We chose $F(x,t)$ and the initial condition in such a way that the exact solution is $u(x,t)=\exp(\cos(x-t))$. The scheme was executed for various values of $N$, specifically $N=24,36,48,72,84$. We used a sixth-order explicit Runge-Kutta method for time integration, with a small time step, and a Final Time of $T=1$.\\
%
%

\begin{figure*}[t!]
    \centering
    \begin{subfigure}[t]{1\textwidth}
        \centering
        \includegraphics[height=7cm]{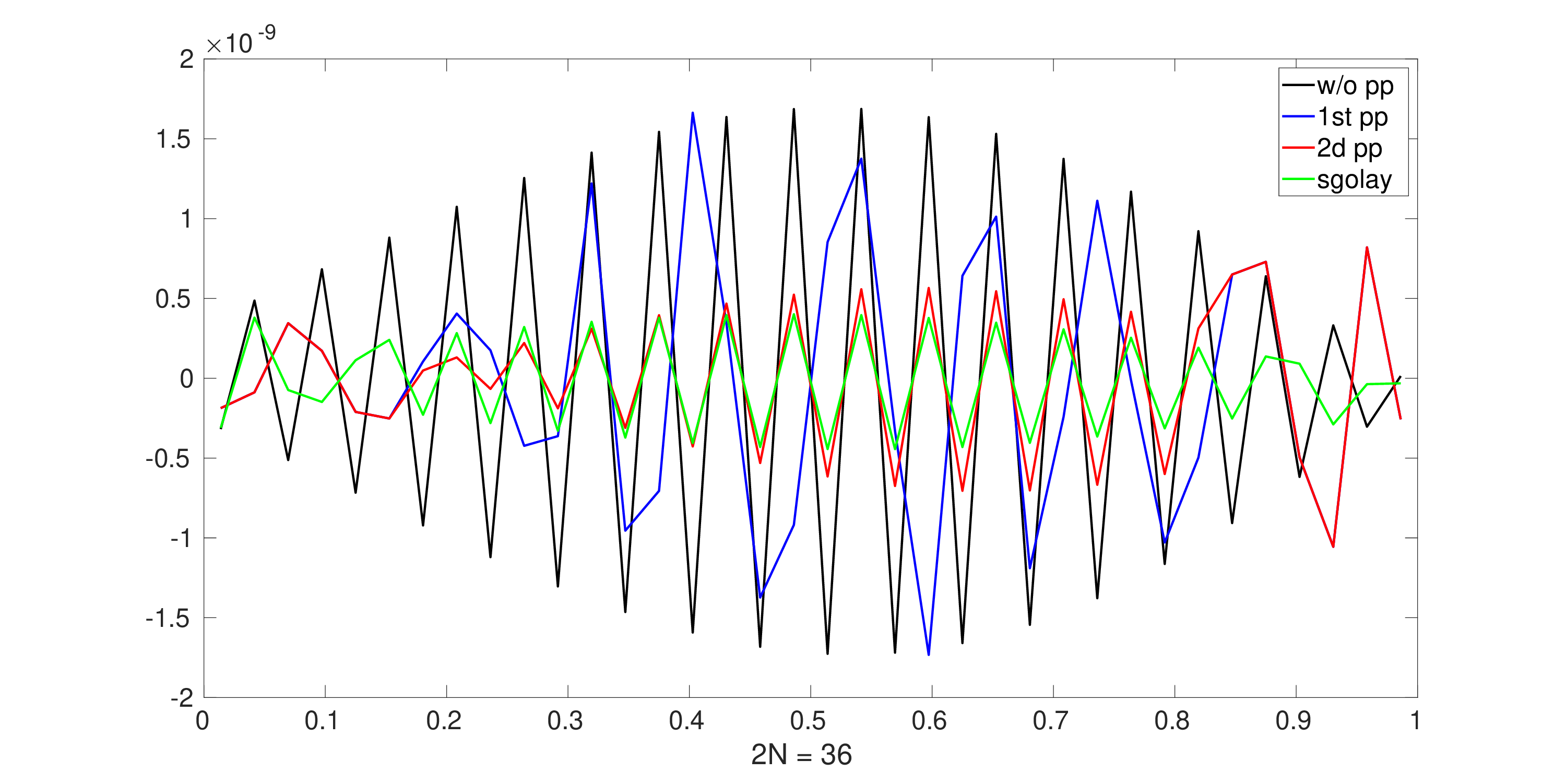}
        \caption{Error plots at Final Time $T=1$ with $N=18$ (36 grid points)}
    \end{subfigure}%
    \quad
    \begin{subfigure}[t]{1\textwidth}
        \centering
        \includegraphics[height=7cm]{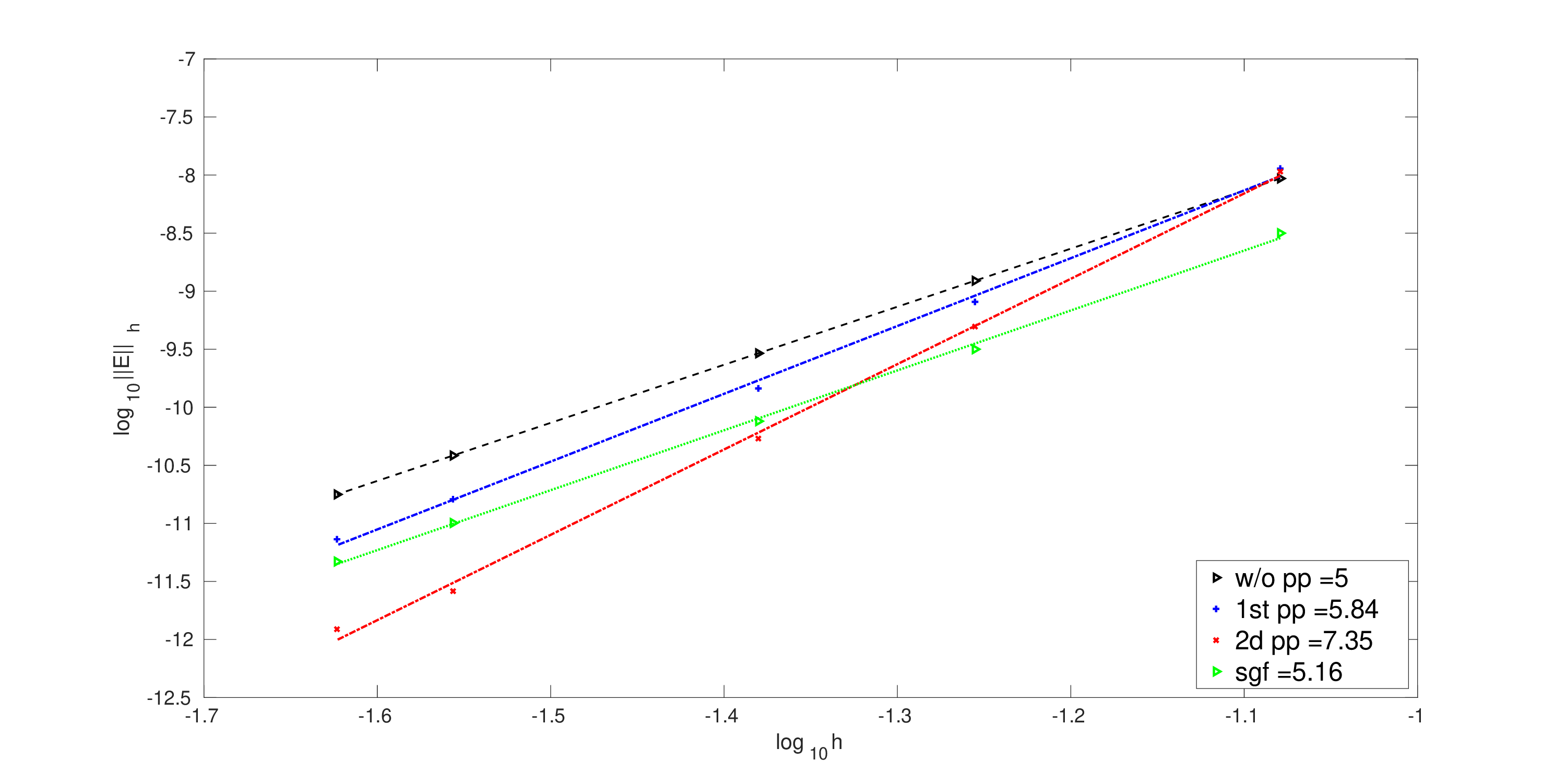}
        \caption{Convergence plots with no post-processing, First and Second Interpolation filters, SG Filter.}
    \end{subfigure}
    \caption{1D Heat Problem - Two Points Block, BFD scheme - Dirichlet BC - $c=-4/13$.}    \label{fig:fig_post_005}
\end{figure*}



The standard SG filter effectively reduces the error but does not affect the order of convergence. Our customized filters are less efficient as far as the error reduction for small $N$, but increase the rate of convergence to seven instead of five, as illustrated in Fig. \ref{fig:fig_post_005}.\\
We note here that further improvements on the issue of filters are left for future research.

\end{document}